\documentclass[11pt,a4paper,twoside]{article}
\usepackage{graphicx}
\usepackage{amsmath}
\usepackage{amsfonts}
\usepackage{amssymb}
\usepackage{enumerate}
\usepackage[hidelinks]{hyperref}
\usepackage{lscape}
\usepackage{longtable}
\usepackage{mathrsfs}
\usepackage{rotating}
\usepackage{color}
\usepackage{url}
\usepackage{rotating}

\usepackage[labelfont=sl,textfont=sl]{subcaption} 
\usepackage[font={small,sl},labelsep=period]{caption} 

\usepackage{titlesec} 
\usepackage[titletoc]{appendix}
\usepackage[ruled,vlined,noline,linesnumbered]{algorithm2e}
\newtheorem{theorem}{Theorem}[section]

\newtheorem{corollary}{Corollary}[section]

\newtheorem{definition}{Definition}[section]

\newtheorem{lemma}{Lemma}[section]

\newtheorem{proposition}{Proposition}[section]
\newtheorem{remark}{Remark}[section]

\newtheorem{assumption}{Assumption}[section]

\newenvironment{proof}[1][Proof]{\textbf{#1.} }{\hfill$\Box$}

\newcommand{\R}{\mathbb{R}}

\newcommand{\N}{\mathbb{N}}

\renewcommand{\P}[1]{\operatorname{\mathbb{P}}\left(#1\right)}
\newcommand{\E}[1]{\operatorname{\mathbb{E}}\left[#1\right]}

\newcommand{\vct}[1]{\boldsymbol{#1}}
\newcommand{\mtx}[1]{\boldsymbol{#1}}

\newcommand{\T}{\mathrm{T}}

\newcommand{\cm}[1]{\operatorname{cm}\left(#1\right)}

\newcommand{\eps}{\epsilon}

\newcommand{\calD}{\mathcal{D}}

\newcommand{\calF}{\mathcal{F}}

\newcommand{\calO}{\mathcal{O}}

\newcommand{\calS}{\mathcal{S}}

\newcommand{\calU}{\mathcal{U}}

\newcommand{\va}{\vct{a}}
\newcommand{\vb}{\vct{b}}

\newcommand{\vd}{\vct{d}}
\newcommand{\ve}{\vct{e}}

\newcommand{\vg}{\vct{g}}

\newcommand{\vp}{\vct{p}}

\newcommand{\vu}{\vct{u}}
\newcommand{\vv}{\vct{v}}

\newcommand{\vx}{\vct{x}}

\newcommand{\vz}{\vct{z}}

\newcommand{\mA}{\mtx{A}}

\newcommand{\mI}{\mtx{I}}

\newcommand{\mP}{\mtx{P}}
\newcommand{\mQ}{\mtx{Q}}
\newcommand{\mR}{\mtx{R}}
\newcommand{\mS}{\mtx{S}}

\newcommand{\mZ}{\mtx{Z}}
\newcommand{\mId}{{\bf I}}

\newcommand{\bigO}{\mathcal{O}} 
\newcommand{\grad}{\nabla}
\newcommand{\gammainc}{\gamma_{\rm inc}}
\newcommand{\gammadec}{\gamma_{\rm dec}}
\newcommand{\flow}{f_{\rm low}}
\newcommand{\kalpha}{\overline{\alpha}_k}
\newcommand{\dmax}{D_{\max}}
\newcommand{\pmax}{P_{\max}}
\newcommand{\mingk}{\tilde{g}_k}

\newcommand{\revised}[1]{\textcolor{black}{#1}}

\title{Direct search based on probabilistic descent in reduced spaces}

\author{Lindon Roberts\thanks{Mathematical Sciences Institute, Building 145, 
Science Road, Australian National University, Canberra ACT 2601, Australia 
(\texttt{lindon.roberts@anu.edu.au}).} 
\and 
Cl\'ement W. Royer\thanks{LAMSADE, CNRS, Universit\'e Paris Dauphine-PSL, 
Place du Mar\'echal de Lattre de Tassigny, 75016 Paris, France 
(\texttt{clement.royer@lamsade.dauphine.fr}). Funding for this author's research 
was partially provided by CNRS INS2I under the grant GASCON and by Agence 
Nationale de la Recherche through program ANR-19-P3IA-0001 (PRAIRIE 3IA 
Institute).}
}

\begin{document}

\maketitle

\begin{abstract}
	Derivative-free algorithms seek the minimum value of a given objective 
	function without using any derivative information. The performance of 
	these methods often worsen as the dimension increases, a phenomenon predicted 
	by their worst-case complexity guarantees. Nevertheless, recent algorithmic 
	proposals have shown that incorporating randomization into otherwise deterministic 
	frameworks could alleviate this effect for direct-search methods.
	\revised{In particular, the best guarantees and practical performance were obtained 
	for direct-search schemes using a random vector uniformly distributed on the 
	sphere and its negative at every iteration. This approach effectively draws 
	directions from a random one-dimensional 
	subspace, yet the properties of such subspaces have not been exploited in direct 
	search, unlike for other derivative-free schemes.}
	\revised{Moreover, existing theory is by design limited to bounded directions, 
	and thus does not fully account for the numerous possibilities for 
    generating random directions (such as drawing from a Gaussian distribution).}
	
	In this paper, we study a generic direct-search algorithm in which the 
	polling directions are defined using random subspaces. Complexity guarantees 
	for such an approach are derived thanks to probabilistic properties 
	related to both the subspaces and the directions used within these 
	subspaces. \revised{Our analysis crucially extends previous deterministic 
    and probabilistic arguments by relaxing the need for directions to be 
    deterministically bounded in norm. As a result, our approach encompasses a wide 
    range of new optimal polling strategies that can be characterized using our subspace 
    and direction properties.} By leveraging results on random subspace 
    embeddings and sketching matrices, we show that better complexity bounds are 
    obtained for randomized instances of our framework. A numerical investigation 
	confirms the benefit of randomization, particularly when done in 
	subspaces, when solving problems of moderately large dimension.
\end{abstract}

\section{Introduction} 
\label{sec:intro}

This paper is concerned with solving the following unconstrained optimization 
problem:
\begin{equation}
\label{eq:pb}
	\min_{\vx \in \R^n} f(\vx),
\end{equation}
where $f: \R^n \rightarrow \R$ is a continuously differentiable function. We 
suppose that the derivative of $f$ is unavailable for algorithmic purposes, 
thereby precluding the use of standard \revised{nonlinear} optimization 
techniques. \revised{Problems of this form commonly arise in complex 
engineering problems where an expensive simulation code is to be calibrated~\cite{CAudet_WHare_2017}. Tackling such problems is the goal of 
derivative-free optimization. Over the past decades, this field has given rise to 
well-established algorithms, that have been successfully applied to
real-world engineering problems}~\cite{CAudet_WHare_2017,
ARConn_KScheinberg_LNVicente_2009b,JLarson_MMenickelly_SMWild_2019}.
\revised{Derivative-free optimization methods can be broadly classified into 
two categories:} model-based algorithms, that construct a model of the 
objective function to guide the optimization process into 
selecting a new point to evaluate, and direct-search methods, that proceed 
by exploring the space along suitably chosen directions.

A recent trend in nonlinear optimization, and in derivative-free optimization 
in particular, is to compare optimization methods based on their 
complexity guarantees.  In the context of derivative-free optimization, we are 
interested in the number of function evaluations required by these methods 
to reach a point $\vx$ such that $\|\nabla f(\vx)\| \le \epsilon$, where 
\revised{$\epsilon >0$} is a given tolerance. Recent results in worst-case 
complexity analyzes have established that this number is of order 
$n^2 \epsilon^{-2}$ 
for standard deterministic derivative-free frameworks such as trust 
region~\cite{RGarmanjani_DJudice_LNVicente_2016} and direct search with 
sufficient decrease~\cite{LNVicente_2013}, as well as $n^2 \epsilon^{-3/2}$ 
for cubic regularization based on finite-difference 
estimates~\cite{CCartis_NIMGould_PhLToint_2012b}. In the case of direct 
search, it was shown that the factor $n^2$ could not be improved 
upon using deterministic algorithms~\cite{MDodangeh_LNVicente_ZZhang_2016}.

Despite this negative result, a number of recent algorithmic proposals 
based on randomized techniques have been shown to reduce the dependency in 
the dimension from $n^2$ to $n$. A direct-search method based on random 
directions was endowed with a complexity bound of order $n \epsilon^{-2}$ 
in a probabilistic sense, provided the directions were chosen uniformly on 
the unit sphere~\cite{SGratton_CWRoyer_LNVicente_ZZhang_2015}. In particular, 
choosing two opposite directions was identified as the best choice from a 
theoretical viewpoint (in that it maximizes the probability of having a 
good direction among two uniformly drawn), that also yielded the best 
results in numerical experiments~\cite{SGratton_CWRoyer_LNVicente_ZZhang_2015,
SGratton_CWRoyer_LNVicente_ZZhang_2019}. Similar results were obtained through 
a different analysis for zeroth-order methods that constructed a gradient 
estimate, typically based on Gaussian 
sampling~\cite{EBergou_EGorbunov_PRichtarik_2020,YuNesterov_2011,
YuNesterov_VSpokoiny_2017}. Very recently, Kozak et 
al.~\cite{DKozak_SBecker_ADoostan_LTenorio_2021} also proposed to use 
directional derivatives sampled within random orthogonal subspaces to 
approximate the gradient. In a related, subsequent work, Kozak et 
al.~\cite{DKozak_CMolinari_LRosasco_LTenorio_SVilla_2021} described a variant 
of this approach using finite difference estimates based on orthogonal 
directions drawn from random subspaces. In all the aforementioned approaches, 
a linear dependence in the problem dimension could be identified, 
representing an improvement over the deterministic setting. 
\revised{Note, however, that the seemingly general analysis proposed in 
the direct-search setting~\cite{SGratton_CWRoyer_LNVicente_ZZhang_2015} did 
not include direction choices such as Gaussian vectors.}

The situation was different in the model-based framework. Although 
randomized frameworks were proposed based on the same reasoning \revised{as} 
direct-search methods, the complexity analysis did not suggest 
any possible improvement in terms of the problem 
dimension~\cite{ASBandeira_KScheinberg_LNVicente_2014,
SGratton_CWRoyer_LNVicente_ZZhang_2018}. However, recent approaches 
focusing on constructing models in random, low-dimensional subspaces 
did manage to achieve such a theoretical 
improvement~\cite{CCartis_LRoberts_2021}
by leveraging random embeddings using for instance Gaussian or random 
orthogonal matrices. Such techniques appeared as promising to develop 
scalable derivative-free optimization techniques, which was done by 
Cartis et al.~\cite{CCartis_TFerguson_LRoberts_2020,CCartis_LRoberts_2021}. Their approach 
builds on a general framework for derivative-based optimization in 
random subspaces~\cite{CCartis_JFowkes_ZShao_2020,CCartis_JFowkes_ZShao_2022,
ZShao_2022}, and adapts it to a derivative-free, model-based setting. However, 
the connection between the subspace arguments in the model-based literature 
and the use of random one-dimensional subspaces in direct search was 
not investigated.

In this paper, we revisit the analysis of direct-search methods based on 
probabilistic properties in order to improve our understanding of their 
behavior. To this end, we propose a framework that relies on search directions 
chosen within subspaces of the variable space. We define probabilistic 
properties on these subspaces and directions, by combining elements 
from the direct-search literature and that of model-based methods based on 
random subspaces. Such properties are then combined to yield complexity 
bounds for our framework, in a way that departs from existing 
analyses~\cite{SGratton_CWRoyer_LNVicente_ZZhang_2015}. In particular, 
our reasoning allows for unbounded directions, and thus encompasses 
popular practical choices such as Gaussian directions. 
\revised{This decoupling of the subspace and direction generation enables us to introduce a suite of new methods that all match the best known 
bound in terms of dependencies on $n$, and find a new interpretation of 
direct search based on opposite, random directions.}
As a result, our framework 
for direct search together with the model-based 
framework~\cite{CCartis_LRoberts_2021} provides a coherent understanding of 
the benefits of randomization for scalable derivative-free optimization of 
nonconvex functions, resolving the disconnect between randomized direct-search 
and randomized model-based techniques that was previously 
observed~\cite{SGratton_CWRoyer_LNVicente_ZZhang_2015,
SGratton_CWRoyer_LNVicente_ZZhang_2018}. Our experiments confirm that a 
randomized subspace approach can be quite beneficial in a direct-search 
setting.

The rest of this document is organized as follows. In 
Section~\ref{sec:ds}, we recall the main features of direct-search 
methods, and the associated complexity guarantees. We then describe a new 
paradigm to generate polling directions based on subspaces in 
Section~\ref{sec:sds}, for which we establish probabilistic 
complexity guarantees.
Numerical experiments for our proposed techniques are given in 
Section~\ref{sec:num}.

\paragraph{Software}
All the algorithms discussed here are available in an open-source Python 
package available on 
Github.\footnote{\url{https://github.com/lindonroberts/directsearch}}

\paragraph{Notation} In what follows, $\|\cdot\|$ will denote the Euclidean 
norm for vectors or the operator 2-norm for matrices. We use $\log(\cdot)$ to 
denote the reciprocal of the exponential function and $\log_{a}(\cdot)$ to 
denote the base-$a$ logarithmic function. The vectors will be denoted 
by bold lowercase letters (e.g. $\vx$) while the matrices will be denoted by 
bold uppercase letters (e.g. $\mS$). Sets will be denoted by cursive uppercase 
letters (e.g. $\calD$). The letters $m,n,r$ will always denote integers 
greater than or equal to $1$. Finally, $\mI_r$ will denote the identity matrix 
in $\R^{r \times r}$.

\section{Direct-search framework and complexity results}
\label{sec:ds}

In this section, we recall the key components of a direct-search method based 
on sufficient decrease. We focus on the properties that guarantee decrease 
in the objective function and, as a result, convergence to approximate 
stationary points. Section~\ref{subsec:ds:det} recalls the fundamentals 
of direct search based on deterministic properties, while 
Section~\ref{subsec:ds:pba} extends the analysis to the case of 
probabilistically descent directions, following the reasoning of Gratton et 
al.~\cite{SGratton_CWRoyer_LNVicente_ZZhang_2015}.

\subsection{Direct search based on deterministic descent}
\label{subsec:ds:det}

Algorithm~\ref{algo:ds} presents a simplified direct-search framework based 
on sufficient decrease, for which complexity results can easily be 
established. At every iteration, the algorithm chooses a set of polling  
directions of fixed cardinality, and evaluates the objective at the 
corresponding points. If one of these trial points satisfies the decrease 
condition~\eqref{eq:ds:suffdec}, this point becomes the new iterate and the 
stepsize parameter is (possibly) increased. Otherwise, the current point does 
not change, and the stepsize is  decreased.
\revised{Such an adaptive behavior of the stepsize parameter is instrumental to 
establishing global convergence of the algorithm, and has roots in convergence 
analyzes based on line-search techniques~\cite{TGKolda_RMLewis_VTorczon_2003}.}

\begin{algorithm}[h!]
	\SetAlgoLined
	\DontPrintSemicolon 
	\BlankLine
	\textbf{Inputs}: $\vx_0 \in \R^n$, $\alpha_{\max} > 0$, 
	$\alpha_0 \in (0,\alpha_{\max}]$, $c>0$, $0 < \gammadec < 1 < \gammainc$, 
	$m \in \N$.\;
	\For{$k=0,1,...$}{
		Compute a polling set $\calD_k \subset \R^n$ of $m$ vectors.\;
		If there exists $\vd_k \in \calD_k$ such that
		\begin{equation}
		\label{eq:ds:suffdec}
			f(\vx_k+\alpha_k \vd_k) < 
			f(\vx_k) - \frac{c}{2} \alpha_k^2 \|\vd_k\|^2,
		\end{equation}
		set $\vx_{k+1}:=\vx_k+\alpha_k \vd_k$ and 
		$\alpha_{k+1}:=\min\{\gammainc \alpha_k,\alpha_{\max}\}$.\;
		Otherwise, set $\vx_{k+1}:=\vx_k$ and $\alpha_{k+1}:=\gammadec \alpha_k$.\;
	}
\caption{Direct-search framework based on sufficient decrease.\label{algo:ds}}
\end{algorithm}

The choice of the polling sets is crucial for obtaining theoretical guarantees on the 
behavior of Algorithm~\ref{algo:ds}. The standard requirements rely on the following 
concept of cosine measure, which we define for an arbitrary dimension $r$ for later 
use in the paper.

\begin{definition}
\label{de:cosmeas}
	Given a set of vectors $\calD \subset \R^r$ and a nonzero vector 
	$\vv \in \R^r$, the \emph{cosine measure of $\calD$ at $\vv$} is 
	defined by
	\begin{equation}
	\label{eq:cosmeasvec}
		\cm{\calD,\vv}:= \max_{\vd \in \calD} \frac{\vd^\T \vv}{\|\vd\|\|\vv\|}.
	\end{equation}
	The \emph{cosine measure of $\calD$} is then given by
	\begin{equation}
	\label{eq:cosmeas}
		\cm{\calD}:=\min_{\substack{\vv \in \R^{r} \\ \|\vv\|\neq 0}} 
		\cm{\calD,\vv}.
	\end{equation}
\end{definition}
Any set $\calD$ such that $\cm{\calD}>0$ is called a positive spanning set (PSS) for 
$\R^r$, as its 
elements span $\R^r$ by nonnegative linear combinations~\cite{CAudet_WHare_2017,ARConn_KScheinberg_LNVicente_2009b}.
Using PSSs leads to complexity results for Algorithm~\ref{algo:ds} under 
the following standard assumptions on the objective function.

\begin{assumption}
\label{as:flow}
    There exists $\flow \in \R$ such that $f(\vx) \ge \flow$ for all $\vx \in \R^n$.
\end{assumption}

\begin{assumption}
\label{as:fC11}
    The function $f$ is continuously differentiable, and its derivative is 
    $L$-Lipschitz continuous with $L>0$.
\end{assumption}

The analysis of Algorithm~\ref{algo:ds} is based on two key arguments. First, one 
can use the sufficient decrease property~\eqref{eq:ds:suffdec} to guarantee that the 
step size converges to zero, \revised{regardless of the cosine measure of the polling 
sets. T}his is the purpose of the following lemma.

\begin{lemma}
\label{le:ds:sumalpha}
    Let Assumption~\ref{as:flow} hold, and consider the step size sequence 
    $\{\alpha_k\}$ produced by Algorithm~\ref{algo:ds}. Suppose that the 
    directions in $\{\calD_k\}_k$ are uniformly bounded in norm. Then, there exists 
    $\beta>0$ that does not depend on $\{\calD_k\}$ such that
    \begin{equation}
    \label{eq:ds:sumalpha}
        \sum_{k=0}^{\infty} \alpha_k^2 \le \beta < \infty.
    \end{equation}
    As a result, $\lim_{k \rightarrow \infty} \alpha_k = 0$. 
\end{lemma}

\revised{A key assumption to derive the result of Lemma~\ref{le:ds:sumalpha} is 
that the polling directions are uniformly bounded in norm. Such a property is 
easy to satisfy in practice (e.g. by normalizing directions), and leads to global convergence~\cite{TGKolda_RMLewis_VTorczon_2003} 
and complexity results~\cite{LNVicente_2013} for deterministic direct search. 
It was also used in analyzing probabilistic direct 
search~\cite{SGratton_CWRoyer_LNVicente_ZZhang_2015}, thereby restricting the 
kind of directions that could be used. Relaxing this requirement in a probabilistic 
fashion leads to a number of complications that we will deal with in Section~\ref{sec:sds}.}


The second main ingredient \revised{in deriving complexity bounds for 
Algorithm~\ref{algo:ds} relates the 
quality of the polling sets to the stepsize. It certifies that the method 
will move to a new point if the stepsize is small compared to the gradient norm 
at the current point.} 

\begin{lemma}
\label{le:ds:unsucc}
	Consider the $k$-th iteration of Algorithm~\ref{algo:ds} under the assumption 
	that\\
	$\cm{\calD_k,-\nabla f(\vx_k)} \ge \kappa \in (0,1)$\revised{and 
	$\|\vd_k\| \le \dmax$ for any $\vd \in \calD_k$. Then, if}
	\begin{equation}
	\label{eq:ds:unsucc}
		\revised{\alpha_k < \frac{2}{(L+c)\dmax}\kappa \|\nabla f(\vx_k)\|,}
	\end{equation}
	\revised{the sufficient decrease condition is satisfied for one direction 
	in $\calD_k$, and thus $\vx_{k+1}\neq \vx_k$.}
\end{lemma}
In practice, the gradient is unknown, thus the assumption 
$\cm{\calD_k,-\nabla f(\vx_k)} \ge \kappa$ is replaced by 
$\cm{\calD_k} \ge \kappa$, which is equivalent to assuming that $\calD_k$ is a 
PSS in $\R^n$.

The updating process on $\{\alpha_k\}_k$ together with the result of 
Lemma~\ref{le:ds:unsucc} guarantees convergence of Algorithm~\ref{algo:ds} 
provided the cosine measure sequence $\{\cm{\calD_k}\}_k$ is uniformly 
bounded below by $\kappa \in (0,1)$.
Under this assumption, it is known~\cite{LNVicente_2013} that 
Algorithm~\ref{algo:ds} reaches an iterate $\vx_k$ 
such that $\|\nabla f(\vx_k)\| \le \epsilon$ using at most 
\begin{equation}
\label{eq:wccdsdeterm}
	\mathcal{O}\left( m\,\kappa^{-2}\,\epsilon^{-2}\right)
\end{equation}
function evaluations.

A classical choice in direct search consists in selecting $\calD_k$ as the 
set of coordinate vectors and their negatives in $\R^n$. In that case, one 
has $m=2n$, $\kappa=\tfrac{1}{\sqrt{n}}$, and the bound becomes 
\begin{equation}
\label{eq:wccdplusdeterm}
	\mathcal{O}\left( n^2\,\epsilon^{-2}\right).
\end{equation}
In a deterministic setting, the dependency in $n^2$ cannot be improved while 
using positive spanning sets without additional 
information~\cite{MDodangeh_LNVicente_ZZhang_2016}.

\subsection{Direct search based on probabilistic descent}
\label{subsec:ds:pba}

As highlighted in the previous section, the use of positive spanning sets is 
instrumental for convergence of classical direct-search methods. However, this 
property also incurs a dependency of $n^2$ in the complexity bounds, due to 
the need to cover an $n$-dimensional space. Gratton et 
al.~\cite{SGratton_CWRoyer_LNVicente_ZZhang_2015} recently established that 
randomly generated direction sets \emph{that do not form a PSS} can still 
provide a good approximation of a particular vector in that space, and that 
this is sufficient to produce good directions. The following property was 
thus introduced.

\begin{definition}
\label{de:ds:pbadesc}
	Given $p \in (0,1]$ and $\kappa \in (0,1)$, the polling set sequence 
	$\{\calD_k\}$ used in Algorithm~\ref{algo:ds} is called 
	$(p,\kappa)$-descent if
	\begin{equation}
	\label{eq:ds:pbadesc}
		\left\{
		\begin{array}{lll}	
			\P{ \cm{\calD_0,-\nabla f(\vx_0)} \ge \kappa} &\ge &p \\
			 & & \\
			\P{ \cm{\calD_k,-\nabla f(\vx_k)} 
			\ge \kappa \middle| \calF_{k-1} } &\ge &p, \qquad \forall k \ge 1,
		\end{array}
		\right.
	\end{equation}
	where $\calF_{k-1}$ is the $\sigma$-algebra generated by the random sets 
	$\calD_0,\dots,\calD_{k-1}$.
\end{definition}

The $(p,\kappa)$-descent property is enough to establish probabilistic 
complexity results for Algorithm~\ref{algo:ds}. Indeed, the result of 
Lemma~\ref{le:ds:sumalpha} holds for every realization of the method provided 
the directions are bounded, while that of Lemma~\ref{le:ds:unsucc} now 
depends on the occurrence of the random event 
$\{\cm{\calD_k,-\nabla f(\vx_k)} \ge \kappa\}$.

Using martingale-type arguments, it is then possible to show that a method 
employing a $(p,\kappa)$-descent sequence converges almost surely to a point 
with zero gradient. In addition, high probability complexity guarantees hold, 
in that the method reaches an iterate satisfying
$\|\nabla f(\vx_k)\|\le \epsilon$ using at most
\[
	\mathcal{O}(m \kappa^{-2}\epsilon^{-2})
\]
function evaluations \emph{with probability at least 
$1-\mathcal{O}(-\exp(C\epsilon^{-2}))$}.
Assuming $\calD_k$ is randomly generated using $m \ge 2$ directions uniformly 
distributed on the unit sphere, one obtains an high-probability evaluation 
complexity bound in 
\begin{equation}
\label{eq:wccdspd}
	\mathcal{O}(n\epsilon^{-2}).
\end{equation}
Using $m=2$ emerged as a good practical alternative, with the use of two 
opposite directions allowing to maximize the probability of having a descent 
set~\cite[Appendix B]{SGratton_CWRoyer_LNVicente_ZZhang_2015}.
Other proposals based on random directions~\cite{YuNesterov_VSpokoiny_2017,
EBergou_EGorbunov_PRichtarik_2020} relied on evaluations along 
random Gaussian opposite directions, leading to a same improvement in the 
complexity bound (note that the probabilistic descent analysis does not 
apply to Gaussian direction since those are not deterministically bounded 
in norm). For such approaches, using evaluations along a one-dimensional 
subspace stood out as an efficient and theoretically sound strategy, but the 
role of the subspace was not further investigated. Taking the subspace nature 
of those directions into account is the key goal of this paper, and the subject 
of the next section.

\section{Probabilistic descent in reduced spaces}
\label{sec:sds}

Building in the framework of Algorithm~\ref{algo:ds}, we propose a method that 
operates in a reduced space by selecting both a subspace of $\R^n$ and a set 
of polling directions within that subspace. This two-step process allows to 
identify deterministic and probabilistic conditions under which the subspace 
(resp. the directions) are of suitable quality.

\subsection{Algorithm and suitable properties}
\label{subsec:sds:algo}

Algorithm~\ref{algo:sds} details our proposed method. At every iteration, we 
first generate directions in $\R^r$ with $r\le n$. We then combine those 
directions with a matrix $\mP_k \in \R^{r \times n}$ to obtain a polling set 
in $\R^n$ \emph{that belongs to an $r$-dimensional subspace of $\R^n$}. This 
property is a key feature of our method, as it allows our method to operate 
in subspaces of lower dimension than that of the original problem: at 
iteration $k$, our polling set is $\{\mP_k^\T \vd | \vd \in \calD_k\}$.

\begin{algorithm}[h!]
	\SetAlgoLined
	\DontPrintSemicolon 
	\BlankLine
	\textbf{Inputs}: $\vx_0 \in \R^n$, $\alpha_{\max} > 0$,  
	$\alpha_0 \in (0,\alpha_{\max}]$, $c>$, $0 < \gammadec < 1 < \gammainc$; 
	$m \in \N$, $r \le n$.\;
	\For{$k=0,1,...$}{
		Compute a matrix $\mP_k \in \R^{r \times n}$.\;
		Compute a set $\calD_k \subset \R^r$ of $m$ vectors.\;
		If there exists $\vd_k \in \calD_k$ such that
		\begin{equation}
		\label{eq:ssuffdec}
			f(\vx_k+\alpha_k \mP^\T_k \vd_k) < f(\vx_k) 
			- \frac{c}{2} \alpha_k^2 \|\mP_k^\T \vd_k\|^2,
		\end{equation}
		set $\vx_{k+1}:=\vx_k+\alpha_k \mP^T_k \vd_k$ and 
		$\alpha_{k+1}:=\min\{\gammainc\alpha_k,\alpha_{\max}\}$.\;
		Otherwise, set $\vx_{k+1}:=\vx_k$ and 
		$\alpha_{k+1}:=\gammadec \alpha_k$.\;
	}
\caption{Direct-search method in reduced spaces.\label{algo:sds}}
\end{algorithm}

To assess the quality of the polling sets, we define separate properties for 
the matrix $\mP_k$ and the set $\calD_k$, starting with the former. The matrix 
$\mP_k$ produces an $r$-dimensional subspace of $\R^n$ in which polling 
directions will be generated. When $r<n$, the use of $\mP_k$ will prevent the 
polling set from providing good approximations to all of $\R^n$. Nevertheless, 
for optimization purposes, we are merely interested in approximating the 
negative gradient (and its norm). Consequently, we require the matrix $\mP_k$ 
to capture a significant portion of gradient information. In addition, 
defining the polling set through application of $\mP_k^\T$ should alter the 
directions in the reduced subspace in a controlled way, which we express 
through bounds on the singular values of the matrix. These considerations lead 
to the following definition, motivated by a similar concept in the model-based 
setting~\cite{CCartis_LRoberts_2021,CCartis_JFowkes_ZShao_2020,
CCartis_JFowkes_ZShao_2022,ZShao_2022}.

\begin{definition}
\label{de:sds:Pk}
	Let $\eta$, $\sigma$ and $\pmax$ be positive quantities.
	For any realization of Algorithm~\ref{algo:sds} and any $k \in \N$, the 
	matrix $\mP_k$ is called \emph{$(\eta,\sigma,\pmax)$-well aligned} for $f$ 
	at $\vx_k$ provided
	\begin{equation}
	\label{eq:sds:Pk:grad}
		\|\mP_k \nabla f(\vx_k)\| \ge \eta \|\nabla f(\vx_k)\|,
	\end{equation}
	\begin{equation}
	\label{eq:sds:Pk:pmax}
		\|\mP_k\| \le \pmax,
	\end{equation}
	\begin{equation}
	\label{eq:sds:Pk:PTd}
		\sigma_{\min}(\mP_k) \ge \sigma,
	\end{equation}
	where $\sigma_{\min}(\cdot)$ denotes the minimum nonzero 
	singular value of the matrix $\mP_k$.
\end{definition}


Conditions~\eqref{eq:sds:Pk:grad}--\eqref{eq:sds:Pk:PTd} are satisfied with 
$\eta=\sigma=\pmax=1$ when $r=n$ and $\mP_k$ is the identity matrix, but may 
not hold when $r< n$, or when $\mP_k$ is a random matrix. For this reason, 
we introduce a probabilistic counterpart to Definition~\ref{de:sds:Pk}.

\begin{definition}
\label{de:sds:pbaPk}
	The sequence $\{\mP_k\}_k$ generated by Algorithm~\ref{algo:sds} is called 
	\emph{$(\eta,\sigma,\pmax,q)$-well-aligned} for $q \in (0,1]$ if
	\begin{equation}
	\label{eq:sds:pbaPk}
		\begin{array}{rll}
			\P{\mP_0\ \mbox{is ($\eta,\sigma,\pmax$)-well aligned}}
			&\ge &q \\
			\forall k \ge 1, \qquad 
			\P{ \mP_k\ \mbox{is ($\eta,\sigma,\pmax$)-well aligned} \ 
			\middle|\ \calF_{k-1}}
			&\ge &q,
		\end{array}
	\end{equation}
	where $\calF_{k-1}$ is the $\sigma$-algebra generated by 
    \revised{$\mP_0,\calD_0,\dots,\mP_{k-1},\calD_{k-1}$}.
\end{definition}

Our requirement on $\calD_k$ is given below, and is similar to that used in 
Section~\ref{subsec:ds:det}.

\begin{definition}
\label{de:sds:Dk}
	Let $\kappa \in (0,1]$ and \revised{$\dmax>1$}. 
	For any realization of Algorithm~\ref{algo:sds} and any index $k \in \N$, 
	the set $\calD_k$ is called $(\kappa,\dmax)$-descent for $f$ and $\mP_k$ at 
	$\vx_k$ if
	\begin{equation}
	\label{eq:sds:Dkcm}
		\cm{\calD_k,-\mP_k \nabla f(x_k)}
		= \max_{\vd \in \calD_k} 
		\frac{-\vd^\T \mP_k \nabla f(\vx_k)}{\|\vd\|\|\mP_k \nabla f(\vx_k)\|} \ge \kappa
	\end{equation}
	and
	\begin{equation}
	\label{eq:sds:Dkdmax}
		\forall \vd \in \calD_k, \quad \dmax^{-1} \le  \|\vd\| \le \dmax.
	\end{equation}
\end{definition}

Examples of sets satisfying Definition~\ref{de:sds:Dk} are positive spanning 
sets in $\R^r$ with unitary elements. As for the properties of $\mP_k$, we provide a 
probabilistic counterpart of Definition~\ref{de:sds:Dk} below.

\begin{definition}
\label{de:sds:pbaDk}
	The sequence \revised{$\{\calD_k\}_k$} generated by Algorithm~\ref{algo:sds} is called 
	\emph{$(\kappa,\dmax,p)$-descent} for $p \in (0,1]$ if
	\begin{equation}
	\label{eq:sds:pbaDk}
		\begin{array}{rll}
			\P{\calD_0\ \mbox{is ($\kappa,\dmax$)-descent} 
			\middle|\ \calF_{-1/2}}
			&\ge &p \\
			\forall k \ge 1, \qquad 
			\P{\calD_k\ \mbox{is ($\kappa,\dmax$)-descent} \ 
			\middle|\ \calF_{k-1/2}}
			&\ge &p,
		\end{array}
	\end{equation}
	where \revised{$\calF_{k-1/2}$} is the $\sigma$-algebra generated by 
	$\mP_0,\calD_0,\dots,\mP_{k-1},\revised{\calD_{k-1}},\mP_k$ and $\calF_{-1/2}$ is 
	the $\sigma$-algebra generated by $\mP_0$.
\end{definition}

Definition~\ref{de:sds:pbaDk} departs from Definition~\ref{de:ds:pbadesc} in 
that it allows for unbounded directions, as long as such directions occur with 
a small probability. This enables for instance the use of Gaussian vectors, a 
distribution that was not immediately covered by the analysis of Gratton et 
al.~\cite{SGratton_CWRoyer_LNVicente_ZZhang_2015}.
In Section~\ref{subsec:sds:ps} we give several possible choices for $\mP_k$ 
and \revised{$\calD_k$} that satisfy these requirements.

\subsection{Complexity analysis}
\label{subsec:sds:wcc}

In this section, we leverage the probabilistic properties defined above to 
derive complexity results for Algorithm~\ref{algo:sds}. Our analysis follows 
a reasoning previously developed for derivative-free methods based on 
probabilistic properties~\cite{SGratton_CWRoyer_LNVicente_ZZhang_2015,
SGratton_CWRoyer_LNVicente_ZZhang_2018}, but involves two properties of this 
form at every iteration, respectively related to $\mP_k$ and $\calD_k$. 
Although these properties can be handled jointly, the analysis still departs 
from existing ones as we allow for directions that are unbounded in norm. At 
the same time, we point out that our approach still relies on exact function 
values, and therefore does not require Lyapunov functions 
similar to those used in stochastic 
optimization~\cite{CPaquette_KScheinberg_2020}.

For the rest of this section, let $\calS$ denote the index set of 
\emph{successful iterations} (i.e.~the $k$ for which $\vx_{k+1} \neq \vx_k$) 
and $\calU$ denote the index set of \emph{unsuccessful iterations} (for which 
$\vx_{k+1}=\vx_k$). \revised{The following lemma describes sufficient 
conditions under which an iteration must be successful, based on our 
properties of interest.}

\begin{lemma}
\label{le:sds:succit}
	Let Assumption~\ref{as:fC11} hold, and consider the $k$-th iteration of 
	a realization of Algorithm~\ref{algo:sds}. Suppose further that 
	$\calD_k$ is $(\kappa,\dmax)$-descent and that $\mP_k$ is 
	$(\eta,\sigma,\pmax)$-well aligned. Finally, suppose that 
	\begin{equation}
	\label{eq:sds:alphabnd}
		\alpha_k < \bar{\alpha}\|\nabla f(\vx_k)\|, \quad \mbox{where} 
		\quad \bar{\alpha}:= \frac{2 \kappa \eta}{(L+c)\pmax^2 \dmax^3}.
	\end{equation}
	Then, the $k$-th iteration is successful.
\end{lemma}

\begin{proof}
	To find a contradiction, suppose that iteration $k$ is unsuccessful.
	Then, by Assumption~\ref{as:fC11}, for all $\vd\in\calD_k$ we have
	\begin{align*}
		-\frac{c}{2}\alpha_k^2 \|\mP_k^T \vd\|^2 
		&\le f(\vx_k+\alpha_k \mP_k^T \vd) - f(\vx_k), \\
		&\le \alpha_k \vd^T \mP_k \grad f(\vx_k) 
		+ \frac{L}{2}\alpha_k^2 \|\mP_k^T \vd\|^2.
	\end{align*}
	Since $\calD_k$ is a $(\kappa,\dmax)$-descent set, there exists 
	$\vd_k\in\calD_k$ such that 
	\begin{align}
		\frac{\vd_k^T (-\mP_k \grad f(\vx_k))}{\|\vd_k\| \: \|\mP_k \grad f(\vx_k)\|} 
		= \cm{\calD_k,-\mP_k \nabla f(\vx_k)} \geq \kappa.
	\end{align}
	Therefore, we obtain
	\begin{eqnarray*}
		-\frac{c}{2}\alpha_k^2 \|\mP_k^T \vd_k\|^2 
		&\le &-\kappa \alpha_k \|\vd_k\|  \|\mP_k \grad f(\vx_k)\| 
		+ \frac{L}{2}\alpha_k^2\|\mP_k^\T \vd_k\|^2 \\
		\kappa \alpha_k \|\vd_k\| \: \|\mP_k \grad f(\vx_k)\|
		&\le & \frac{L+c}{2} \alpha_k^2\|\mP_k^\T \vd_k\|^2 \\
		\kappa \|\vd_k\| \: \|\mP_k \grad f(\vx_k)\|
		&\le & \frac{L+c}{2} \alpha_k\|\mP_k^\T \vd_k\|^2 \\
	\end{eqnarray*}
	Using now the property~\eqref{eq:sds:Pk:grad} on $\mP_k$ together with 
	the bound~\eqref{eq:sds:Dkdmax} on $\|\vd_k\|$ leads to
	\begin{equation*}
		\kappa \|\vd_k\| \: \|\mP_k \grad f(\vx_k)\|
		\ge \kappa \eta \dmax^{-1} \|\grad f(\vx_k)\|
	\end{equation*}
	as well as
	\begin{equation*}
		\|\mP_k^\T \vd_k\|^2 \le \|\mP_k^\T\|^2 \|\vd_k\|^2 
		\le \pmax^2 \dmax^2.
	\end{equation*}
	Putting everything together, we arrive at
	\begin{equation*}
		\kappa \eta \dmax^{-1} \|\grad f(\vx_k)\| 
		\le \frac{L+c}{2} \pmax^2 \dmax^2 \alpha_k \quad 
		\Leftrightarrow \quad \alpha_k \ge \bar{\alpha}\|\grad f(\vx_k)\|,
	\end{equation*}
	and this contradicts~\eqref{eq:sds:alphabnd}.
\end{proof}

We now introduce the following indicator variables:
\begin{subequations} 
\label{eq:sds:indic}
	\begin{align}
		\label{eq:sds:indiczk}
		Z_k &:= \mathbf{1}\left( \calD_k\ \mbox{($\kappa,\dmax$)-descent and }
		\ \mP_k\  \mbox{($\eta,\sigma,\pmax$)-well aligned}
		\right), \\
		\label{eq:sds:indicvk}
		V_k(\alpha) &:= \mathbf{1}\left( \alpha_k < \alpha \right) 
		\quad \forall \alpha>0, \\
		\label{eq:sds:indicwk}
		W_k &:= \mathbf{1}\left( k \in \calS\right),
	\end{align}
\end{subequations}
whose realizations will be denoted by $z_k,v_k(\alpha),w_k$, respectively. Our 
goal is to bound the sum of $z_k$ by a quantity involving the gradient 
norm. To this end, we study careful combinations of the indicator variables 
above, for which we can provide bounds that are independent of the iteration: 
this technique has proven useful in establishing complexity guarantees for 
derivative-free optimization schemes with probabilistic 
components~\cite{CCartis_LRoberts_2021,CCartis_KScheinberg_2018,
CPaquette_KScheinberg_2020,SGratton_CWRoyer_LNVicente_ZZhang_2015,
SGratton_CWRoyer_LNVicente_ZZhang_2018}.

Our first result bounds the sum of squared stepsizes for a certain 
subset of successful iterations. Unlike in analyses of direct search 
based solely on descent 
properties~\cite{SGratton_CWRoyer_LNVicente_ZZhang_2015}, where such a 
bound can be obtained for all iterations, here it only holds for those 
successful iterations that use both well-aligned subspace matrices and 
descent directions.

\begin{lemma}
\label{le:sds:sumzalpha}
	Let Assumption~\ref{as:flow} hold. For any realization of 
	Algorithm~\ref{algo:sds},
	\begin{equation}
	\label{eq:sds:sumzalpha}
		\sum_{k=0}^{\infty} z_k w_k \alpha_k^2 \le 
		\beta:=\frac{2 \dmax^2(f(\vx_0)-\flow)}{c \sigma^2}.
	\end{equation}
\end{lemma}

\begin{proof}
	It suffices to consider iterations for which $z_k w_k=1$, i.e. 
	successful iterations for which $\calD_k$ is $(\kappa,\dmax)$-descent 
	and $\mP_k$ is $(\eta,\sigma,\pmax)$-well aligned. For such an iteration, 
	there exists $\vd_k \in \calD_k$ such that 
	$\vx_{k+1}=\vx_k+\alpha_k \mP_k^\T \vd_k$, and
	\begin{eqnarray*}
		f(\vx_k) - f(\vx_{k+1}) 
		\ge \frac{c}{2}\alpha_k^2 \|\mP_k^\T \vd_k\|^2 
		\ge \frac{c}{2}\alpha_k^2 \sigma^2 \|\vd_k\|^2 
		\revised{\ge \frac{c}{2}\sigma^2 \dmax^{-2} \alpha_k^2.}
	\end{eqnarray*}
	On the other hand, Assumption~\ref{as:flow} guarantees that
	\begin{equation*}
		f(\vx_0) - \flow 
		\ge \sum_{k=0}^{\infty} f(\vx_k)-f(\vx_{k+1}) 
		\ge \sum_{k=0}^{\infty} z_k w_k\,(f(\vx_k) - f(\vx_{k+1})). 
	\end{equation*}
	Thus, we obtain
	\[
		\sum_{k=0}^{\infty} z_k  w_k \alpha_k^2  
		\le \frac{2 \dmax^2(f(\vx_0)-\flow)}{c \sigma^2},
	\]
	proving the desired result.
\end{proof}

The next two results are obtained by carefully examining the behavior 
of the stepsize sequence. This is another notable departure from 
the analysis of direct search based on probabilistic 
descent~\cite{SGratton_CWRoyer_LNVicente_ZZhang_2015}, that is 
due to our two probabilistic properties. We note that similar 
results have been derived in the context of randomized model-based 
methods~\cite{CCartis_LRoberts_2021}.

\begin{lemma}
\label{le:sds:bndwv}
	Let Assumption~\ref{as:flow} hold. Consider a realization of 
	Algorithm~\ref{algo:sds} and an index $k$ such that 
	$\min_{0 \le \ell \le k-1} \|\nabla f(\vx_{\ell})\| > 0$. 
	Then,
	\begin{equation}
	\label{eq:sds:bndwv}
		\sum_{j=0}^{k-1} v_j\left(\kalpha\right) w_j 
		\; \le \;
		\mu \sum_{j=0}^{k-1} 
		v_j\left(\frac{\gammainc}{\gammadec}\kalpha\right)
		(1-w_j)
	\end{equation}
	with $\mu:=\log_{\gammainc}(\gammadec^{-1})$.	
	\begin{equation}
	\label{eq:sds:kalpha}
		\kalpha:=\frac{\gammadec}{\gammainc}\min\left\{
		\gammainc^{-1} \alpha_0,\bar{\alpha} 
		\min_{0 \le \ell \le k-1} \|\nabla f(\vx_{\ell})\| \right\},
	\end{equation}
	and $\bar{\alpha}$ is defined as in Lemma~\ref{le:sds:succit}.
\end{lemma}

\begin{proof}
	If $v_j\left(\kalpha\right) w_j=0$ for every $j \leq k-1$, the 
	bound clearly holds. For the rest of the proof, we thus 
	assume that there exist at least one index $j \in \{0,\dots,k-1\}$ 
	such that $v_j\left(\kalpha\right) w_j=1$. Since 
	$\tfrac{\gammainc}{\gammadec} > 1$, we also have 
	$v_j\left(\tfrac{\gammainc}{\gammadec}\kalpha\right)w_j=1$, 
	thus there also exists at least one index satisfying this property.
	
	Consider a sequence of iterates $j_1,\dots,j_2$ such that 
	$v_j\left(\tfrac{\gammainc}{\gammadec}\kalpha\right) =1$ for 
	every \revised{$j \in \{ j_1,\dots,j_2\}$}, with
	$v_{j_2}\left(\kalpha\right)w_{j_2}=1$ and either $j_1=0$ or 
	$v_{j_1-1}\left(\tfrac{\gammainc}{\gammadec}\kalpha\right)=0$ 
	(note that such a sequence necessarily exists by assumption). 
	Using the updating rules on the stepsize together with 
	the bound 
	$\kalpha < \gammainc^{-1}\alpha_0 \leq \gammainc^{-1}\alpha_{\max}$, 
	we obtain for any $j=j_1,\dots,j_2$ that
	\begin{equation}
	\label{eq:sds:alphaupvw}
		\left\{
			\begin{array}{lll}
				\alpha_{j+1} 
				&= \min\{\gammainc \alpha_j,\alpha_{\max}\} 
				= \gammainc \alpha_j, 
				&\mbox{if\ } v_j(\kalpha) w_j=1, \\
				\alpha_{j+1}
				&\ge \alpha_j, 
				&\mbox{if\ } (1-v_j(\kalpha))
				v_j\left(\tfrac{\gammainc}{\gammadec}\kalpha\right) w_j=1, \\
				\alpha_{j+1} 
				&= \gammadec \alpha_j, 
				&\mbox{if\ } v_j\left(\tfrac{\gammainc}{\gammadec}\kalpha\right)(1-w_j)=1.
			\end{array}
		\right.
	\end{equation}	
	Note that the rules~\eqref{eq:sds:alphaupvw} cover all possible 
	cases. Applying rule~\eqref{eq:sds:alphaupvw} iteratively gives:
	\begin{eqnarray*}
	    \alpha_{j_2+1} 
	    &\ge &\gammainc^{\sum_{j=j_1}^{j_2} v_j(\kalpha) w_j} 
	    \gammadec^{\sum_{j=j_1}^{j_2} 
	    v_j\left(\tfrac{\gammainc}{\gammadec}\kalpha\right)(1-w_j)} \alpha_{j_1}.
	\end{eqnarray*}
	Moreover, since 
	$v_{j_2}\left(\kalpha\right) w_{j_2}=1$, we have
	\[
	    \alpha_{j_2+1} = \gammainc \alpha_{j_2} \le \gammainc \kalpha,
	\]
	leading to
	\begin{equation}
	\label{eq:vwboundsubseq}
	    \frac{\gammainc\kalpha}{\alpha_{j_1}} 
	    \ge \gammainc^{\sum_{j=j_1}^{j_2} v_j(\kalpha) w_j} 
	    \gammadec^{\sum_{j=j_1}^{j_2} 
	    v_j\left(\tfrac{\gammainc}{\gammadec}\kalpha\right)(1-w_j)}.
	\end{equation}
	
	To bound the ratio $\tfrac{\gammainc\kalpha}{\alpha_{j_1}}$, we consider 
	two cases. If $j_1=0$, then $\gammadec\kalpha \le \alpha_0$ by definition of 
	$\kalpha$, and thus $\gammadec\frac{\kalpha}{\alpha_{j_1}} \le 1$. Otherwise, 
	if $j_1>0$, we have by definition of $j_1$ that 
	$v_{j_1-1}\left(\tfrac{\gammainc}{\gammadec}\kalpha\right)=0$, that is, 
	$\alpha_{j_1-1} \ge \frac{\gammainc}{\gammadec}\kalpha$. Per the updating rules 
	on the step size, this implies
	\[
		\alpha_{j_1} \ge \gammadec \alpha_{j_1-1} \ge \gammainc \kalpha,
	\]
	hence $\gammainc\frac{\kalpha}{\alpha_{j_1}} \le 1$ also holds in this case. 
	Plugging this bound into\revised{~\eqref{eq:vwboundsubseq}} yields
	\[
		1 \ge \gammainc^{\sum_{j=j_1}^{j_2} v_j(\kalpha) w_j} 
	    \gammadec^{\sum_{j=j_1}^{j_2} 
	    v_j\left(\tfrac{\gammainc}{\gammadec}\kalpha\right)(1-w_j)}.
	\]	
	
	Taking logarithms, we obtain that
	\begin{equation*}
	    0 \ge \log(\gammainc)\sum_{j=j_1}^{j_2} v_j(\kalpha) w_j 
		+ \log(\gammadec)\sum_{j=j_1}^{j_2}  
		v_j\left(\tfrac{\gammainc}{\gammadec}\kalpha\right)(1-w_j),
	\end{equation*}
	which after rearranging becomes
	\begin{equation*}
		\sum_{j=j_1}^{j_2} v_j(\kalpha) w_j
		\le  \log_{\gammainc}(\gammadec^{-1})
	    \sum_{j=j_1}^{j_2}  v_j\left(\tfrac{\gammainc}{\gammadec}\kalpha\right)
	    (1-w_j)
	    = \mu \sum_{j=j_1}^{j_2}  v_j\left(\tfrac{\gammainc}{\gammadec}\kalpha\right)
	    (1-w_j).
	\end{equation*}
	To conclude, we simply observe that
	\[
		\left\{ j\ \middle|\ v_j(\kalpha) w_j=1 \right\} 
		\; \subset \; 
		\left\{ j\ \middle|\ v_j\left(\tfrac{\gammainc}{\gammadec}\kalpha\right)=1 \right\},
	\]
	\revised{or} equivalently, that every $j\in\{0,\ldots,k-1\}$ for which $v_j(\kalpha) w_j=1$ 
	is in such a subsequence $j_1,\ldots,j_2$. As a result, we obtain
	\[
		\sum_{j=0}^{k-1} v_j(\kalpha) w_j 
		\; \le \;
		\mu \sum_{j=0}^{k-1} v_j\left(\tfrac{\gammainc}{\gammadec}\kalpha\right)
	    (1-w_j).
	\]	
\end{proof}

\begin{lemma}
\label{le:sds:bnd1w1v}
	Let Assumption~\ref{as:flow} hold. Consider a realization of 
	Algorithm~\ref{algo:sds} and an index $k$ such that 
	$\min_{0 \le \ell \le k-1} \|\nabla f(\vx_{\ell})\| >0$. Then,
	\begin{equation}
	\label{eq:sds:bnd1w1v}
		\sum_{j=0}^{k-1} 
		\left(1-v_j\left(\revised{\kalpha}\right)\right)(1-w_j)
		\; \le \;
		\frac{1}{\mu}\sum_{j=0}^{k-1} \left(1-v_j\left(
		\frac{\gammadec}{\gammainc}\revised{\kalpha}\right)\right)
		w_j
		+ \log_{\gammadec^{-1}}\left(\frac{\alpha_0}{\gammadec\kalpha} \right),
	\end{equation}
	where $\mu$ and $\kalpha$ are defined as in 
	Lemma~\ref{le:sds:bndwv}.
\end{lemma}

\begin{proof}
	The proof follows the template of that of Lemma~\ref{le:sds:bndwv}.
	The bound trivially holds if 
	$\left(1-v_j\left(\revised{\kalpha}\right)\right)(1-w_{j}) =0$ 
	for every $j \leq k-1$. Therefore, we suppose that there exists at least one 
	index $j \in \{0,\dots,k-1\}$ such that 
	$(1-v_j(\kalpha))(1-w_j)=1$. Since $\tfrac{\gammadec}{\gammainc} < 1$, we 
	also have 
	$\left(1-v_j\left(
	\tfrac{\gammadec}{\gammainc}\kalpha\right)\right)(1-w_j)=1$, thus 
	there also exists at least one index satisfying this property.
	
	Consider now a sequence of iterates $j_1,\dots,j_2$ such that 
	$\left(1-v_j\left(
	\tfrac{\gammadec}{\gammainc}\kalpha\right)\right)=1$ 
	for every $j=j_1,\dots,j_2$, with 
	$\left(1-v_{j_2}\left(\kalpha\right)\right)(1-w_{j_2})=1$ and 
	either $j_1=0$ or $v_{j_1-1}\left(
	\tfrac{\gammadec}{\gammainc}\kalpha\right)=1$: such a sequence
	necessarily exists by assumption.
	From the updating rules on the stepsize, we 
	obtain the following possible cases:
	\begin{equation}
	\label{eq:sds:alphaup1v1w}
		\left\{
			\begin{array}{lll}
				\alpha_{j+1} 
				&= \gammadec \alpha_j ,
				&\mbox{if\ } (1-v_j(\kalpha)) (1-w_j)=1, \\
				\alpha_{j+1}
				&\le \alpha_j, 
				&\mbox{if\ } v_j(\kalpha)
				\left(1-v_j\left(\tfrac{\gammadec}{\gammainc}\kalpha\right)\right) 
				(1-w_j)=1, \\
				\alpha_{j+1} 
				&\le \gammainc \alpha_j, 
				&\mbox{if\ } 
				\left(1-v_j\left(\tfrac{\gammadec}{\gammainc}\kalpha\right)\right)
				w_j=1,
			\end{array}
		\right.
	\end{equation}		
	for any $j=j_1,\dots,j_2$. By applying rule~\eqref{eq:sds:alphaup1v1w} 
	iteratively, we  thus obtain
	\[
		\alpha_{j_2+1} \le 
		\gammadec^{\sum_{j=j_1}^{j_2} (1-v_j(\kalpha)) (1-w_j)}
		\gammainc^{\sum_{j=j_1}^{j_2} 
		\left(1-v_j\left(\tfrac{\gammadec}{\gammainc}\kalpha\right)\right)w_j}
		\alpha_{j_1}.
	\]
	In addition, using that $\left(1-v_{j_2}\left(\kalpha\right)\right)(1-w_{j_2})=1$, 
	we also have
	\[
		\alpha_{j_2+1} = \gammadec \alpha_{j_2} 
		\ge \gammadec \kalpha,
	\]
	thus
	\begin{equation}
	\label{eq:sds:1v1wboundsubseq}
		\gammadec\frac{\kalpha}{\alpha_{j_1}}
		\le 
		\gammadec^{\sum_{j=j_1}^{j_2} (1-v_j(\kalpha)) (1-w_j)}
		\gammainc^{\sum_{j=j_1}^{j_2} 
		\left(1-v_j\left(\tfrac{\gammadec}{\gammainc}\kalpha\right)\right)w_j}.
	\end{equation}

	Taking the logarithm, we get
	\begin{equation*}
	   	\log\left(\gammadec\frac{\kalpha}{\alpha_{j_1}}\right) \le 
		-\log(\gammadec^{-1}) \sum_{j=j_1}^{j_2} (1-v_j(\kalpha)) (1-w_j) 
		+ \log(\gammainc) \sum_{j=j_1}^{j_2} 
		\left(1-v_j\left(\tfrac{\gammadec}{\gammainc}\kalpha\right)\right)w_j, 
	\end{equation*}
	which after re-arranging gives
	\begin{eqnarray*}
		\sum_{j=j_1}^{j_2} (1-v_j(\kalpha)) (1-w_j) 
		&\le 
		&\log_{\gammadec^{-1}}(\gammainc)\sum_{j=j_1}^{j_2} 
		\left(1-v_j\left(\tfrac{\gammadec}{\gammainc}\kalpha\right)\right)w_j 
		- \log_{\gammadec^{-1}}\left(\gammadec\frac{\kalpha}{\alpha_{j_1}}\right) \\
		&= &\frac{1}{\mu}\sum_{j=j_1}^{j_2} 
		\left(1-v_j\left(\tfrac{\gammadec}{\gammainc}\kalpha\right)\right)w_j
		+ \log_{\gammadec^{-1}}\left(\frac{\alpha_{j_1}}{\gammadec\kalpha}\right) 
	\end{eqnarray*}
	by definition of $\mu$. 
	
	To bound the last term, we consider two cases.
	If $j_1=0$, then $\frac{\alpha_{j_1}}{\gammadec\kalpha} 
	= \frac{\alpha_0}{\gammadec\kalpha}$. Otherwise, 
	if $j_1>0$, we have by definition of $j_1$ that 
	$v_{j_1-1}\left(\tfrac{\gammadec}{\gammainc}\kalpha\right)=1$, that is, 
	$\alpha_{j_1-1} < \frac{\gammadec}{\gammainc}\kalpha$. Per the updating rules 
	on the step size, this implies
	\[
		\alpha_{j_1} \le \gammainc \alpha_{j_1-1} < \gammadec \kalpha,
	\]
	hence $\frac{\alpha_{j_1}}{\gammadec\kalpha} \le 1$ and thus
	$\log_{\gammadec^{-1}}\left(\frac{\alpha_{j_1}}{\gammadec\kalpha}\right) < 0$.
	Overall, we thus obtain that
	\[
		\sum_{j=j_1}^{j_2} (1-v_j(\kalpha)) (1-w_j) 
		\le 
		\left\{
		\begin{array}{ll}		
			\frac{1}{\mu}\sum_{j=j_1}^{j_2} 
			\left(1-v_j\left(\tfrac{\gammadec}{\gammainc}\kalpha\right)\right)w_j 
			+ \log_{\gammadec^{-1}}\left(\frac{\alpha_0}{\gammadec\kalpha} \right)
			&\mbox{if\ $j_1=0$} \\
			\frac{1}{\mu}\sum_{j=j_1}^{j_2} 
			\left(1-v_j\left(\tfrac{\gammadec}{\gammainc}\kalpha\right)\right)w_j
			&\mbox{otherwise.}
		\end{array}
		\right.
	\]	
	Finally, we note that
	\[
		\left\{ j\ \middle|\ (1-v_j(\kalpha)) (1-w_j)=1 \right\} 
		\; \subset \; 
		\left\{ j\ \middle|\ 1-v_j\left(\tfrac{\gammadec}{\gammainc}\kalpha\right)=1 \right\},
	\]
	on equivalently, that every $j\in\{0,\ldots,k-1\}$ for which $(1-v_j(\kalpha)) (1-w_j)=1$ 
	is in such a subsequence $j_1,\ldots,j_2$.
	This allows us to conclude
	\[
		\sum_{j=0}^{k-1} (1-v_j(\kalpha)) (1-w_j) 
		\le \frac{1}{\mu}\sum_{j=0}^{k-1} 
		\left(1-v_j\left(\tfrac{\gammadec}{\gammainc}\kalpha\right)\right)w_j 
		+ \log_{\gammadec^{-1}}\left(\frac{\alpha_0}{\gammadec\kalpha} \right).
	\]	
\end{proof}

The results of Lemmas~\ref{le:sds:bndwv} and~\ref{le:sds:bnd1w1v} are 
sufficient to obtain a bound on the number of iterations for which the 
directions are generated from both a descent set and a well-aligned subspace 
matrix.

\begin{proposition}
\label{pr:sds:bndz}
	For any realization of Algorithm~\ref{algo:sds} and any positive integer 
	$k$,
	\begin{equation}
	\label{eq:sds:sumindic}
		\sum_{j=0}^{k-1} z_{j} \le 
		\frac{(1-p_0)C}
	    {\min\{\gammainc^{-2}\alpha_0^2,
	    \bar{\alpha}^2 [\mingk]^2\}}
	    + (1-p_0)\log_{\gammadec^{-1}}\left(\frac{\gammainc\alpha_0}{\gammadec^2
	    \min\{\gammainc^{-1}\alpha_0,
	    \bar{\alpha} \mingk\}} \right)
	    + p_0 k,
	\end{equation}
	where $\mingk = \min_{0 \le j \le k-1} \|\nabla f(\vx_j)\|$,
	\begin{equation}
	\label{eq:sds:constant}
		C := 
		\frac{(\mu\gammadec^2+\gammainc^2)\gammainc^2 \beta}{ \mu\gammadec^4},
	\end{equation}
	and 
	\begin{equation}
	\label{eq:sds:p0}
	    p_0 := \max\left\{
    	\frac{\ln(\gammadec)}{\ln(\gammainc^{-1}\gammadec)},
    	\frac{\ln(\gammainc)}{\ln(\gammainc\gammadec^{-1})} \right\} = 
    	\max\left\{\frac{1}{1+\mu},\frac{\mu}{1+\mu}\right\}.
	\end{equation}
\end{proposition}

\begin{proof}
	For any $j=0,\dots,k-1$, we have 
	\begin{equation}
	\label{eq:sds:zdecomp}
		z_j = z_j v_j(\kalpha) w_j + z_j (1-v_j(\kalpha))(1-w_j)
		+ z_j (1-v_j(\kalpha))w_j.
	\end{equation}
	Indeed, the bound clearly holds when $z_j=0$. Assuming $z_j=1$, 
	it also holds when $v_j(\kalpha) w_j=1$ or $(1-v_j(\kalpha))(1-w_j)=1$, while 
	Lemma~\ref{le:sds:succit} implies that we cannot have $z_j v_j(\kalpha) (1-w_j)=1$. This last property will be instrumental in the proof 
	of our result.
	
	Summing~\eqref{eq:sds:zdecomp} over all $j=0,\dots,k-1$, we obtain:
	\begin{equation}
	\label{eq:sds:sumzdecomp}
		\sum_{j=0}^{k-1} z_j 
		= \sum_{j=0}^{k-1} z_j (1-v_j(\kalpha))w_j
		+ \sum_{j=0}^{k-1} z_j v_j(\kalpha) w_j 
		+ \sum_{j=0}^{k-1} z_j (1-v_j(\kalpha)) (1-w_j).
	\end{equation}
	We will provide separate bounds on the three sums on the right-hand side 
	of~\eqref{eq:sds:sumzdecomp}. 
	
	Consider first an index such that $z_j\,(1-v_j(\kalpha))\,w_j=1$. By 
	definition of $v_j(\kalpha)$, we obtain from 
	$1-v_j(\kalpha)=1$ that
	\[
		1 
		\le \left(\frac{\alpha_j}{\kalpha}\right)^2
		= \frac{\gammainc^2\alpha_j^2}{\gammadec^2\min\{\gammainc^{-2}\alpha_0^2,
		\bar{\alpha}^2 \min_{0 \le \ell \le k-1}\|\nabla f(\vx_{\ell})\|^2\}}.
	\]
	Thus,
	\begin{eqnarray*}
		z_j\,(1-v_j(\kalpha))\,w_j
		\le \frac{z_j w_j \gammainc^2\alpha_j^2}{\gammadec^2\min\{
		\gammainc^{-2}\alpha_0^2,
		\bar{\alpha}^2 [\mingk]^2\}} \\
	\end{eqnarray*}
	and summing over all indices up to $k-1$ gives
	\begin{equation}
	\label{eq:sds:sumz1vw}
		\sum_{j=0}^{k-1} z_j\,(1-v_j(\kalpha))\,w_j 
		\le \frac{\gammainc^2}{\gammadec^2\min\{
		\gammainc^{-2}\alpha_0^2,
		\bar{\alpha}^2 [\mingk]^2\}}
		\sum_{j=0}^{k-1} z_j w_j \alpha_j^2
		\le \frac{\gammainc^2 \beta}{\gammadec^2\min\{
		\gammainc^{-2}\alpha_0^2,
		\bar{\alpha}^2 [\mingk]^2\}},
	\end{equation}
	where the last inequality follows from Lemma~\ref{le:sds:sumzalpha}.

	We now bound the second term on the right-hand side 
	of~\eqref{eq:sds:sumzdecomp}. Thanks to Lemma~\ref{le:sds:bndwv}, we 
	have
	\begin{eqnarray*}
		\sum_{j=0}^{k-1} z_j v_j(\kalpha) w_j 
		&\le &\sum_{j=0}^{k-1} v_j(\kalpha) w_j \\
		&\le 
		&\mu \sum_{j=0}^{k-1} v_j\left(\frac{\gammainc}{\gammadec}
		\kalpha\right) (1-w_j) \\
		&= 
		&\mu\left[ 
		\sum_{j=0}^{k-1} z_j v_j\left(\frac{\gammainc}{\gammadec}
		\kalpha\right) (1-w_j) 
		+ 
		\sum_{j=0}^{k-1} (1-z_j) v_j\left(\frac{\gammainc}{\gammadec}
		\kalpha\right) (1-w_j) \right] \\
		&\le 
		&\mu\left[
		\sum_{j=0}^{k-1} z_j v_j\left(\frac{\gammainc}{\gammadec}
		\kalpha\right) (1-w_j) 
		+ 
		\sum_{j=0}^{k-1} (1-z_j) (1-w_j) \right] \\
	\end{eqnarray*}
	Now, for any $j=0,\dots,k-1$, we have
	\[
		z_j v_j\left(\frac{\gammainc}{\gammadec}
		\kalpha\right) (1-w_j)
		\le z_j v_j(\bar{\alpha} \mingk) (1-w_j) = 0,
	\]
	hence the first sum in the above expression is always zero. 
	Thus,
	\begin{equation}
	\label{eq:sds:sumzvw}
	    \sum_{j=0}^{k-1} z_j v_j(\kalpha) w_j 
	    \le 
	    \mu \sum_{j=0}^{k-1} (1-z_j)(1-w_j).
	\end{equation}
	
	We finally consider the last sum in the right-hand side 
	of~\eqref{eq:sds:sumzdecomp}. Applying Lemma~\ref{le:sds:bnd1w1v}, 
	we get
	\begin{eqnarray*}
		\sum_{j=0}^{k-1} z_j (1-v_j(\kalpha)) (1-w_j) 
		&\le 
		&\sum_{j=0}^{k-1} (1-v_j(\kalpha)) (1-w_j) \\
		&\le 
		&\frac{1}{\mu}\sum_{j=0}^{k-1} \left(1-v_j\left(
		\frac{\gammadec}{\gammainc}\kalpha\right)\right)
		w_j + \log_{\gammadec^{-1}}\left(\frac{\alpha_0}{\gammadec\kalpha} \right)\\
		&= 
		&\frac{1}{\mu}\sum_{j=0}^{k-1} z_j \left(1-v_j\left(
		\frac{\gammadec}{\gammainc}\kalpha\right)\right)w_j \\
		&
		&+\frac{1}{\mu}\sum_{j=0}^{k-1} (1-z_j) \left(1-v_j\left(
		\frac{\gammadec}{\gammainc}\kalpha\right)\right)w_j 
		+ \log_{\gammadec^{-1}}\left(\frac{\alpha_0}{\gammadec\kalpha} \right)\\
		&\le 
		&\frac{1}{\mu}\sum_{j=0}^{k-1} z_j \left(1-v_j\left(
		\frac{\gammadec}{\gammainc}\kalpha\right)\right)w_j
		\\
		&
		&+\frac{1}{\mu}\sum_{j=0}^{k-1} (1-z_j) w_j
		+ \log_{\gammadec^{-1}}\left(\frac{\alpha_0}{\gammadec\kalpha} \right).
	\end{eqnarray*}
	
	By the same reasoning used to obtain~\eqref{eq:sds:sumz1vw}, we can 
	bound the first sum in the last expression as follows:
	\[
		\sum_{j=0}^{k-1} z_j \left(1-v_j\left(
		\frac{\gammadec}{\gammainc}\kalpha\right)\right)w_j 
		\le 
		\frac{\gammainc^4 \beta}{\gammadec^4\min\{
		\gammainc^{-2}\alpha_0^2,
		\bar{\alpha}^2 \|\tilde{\vg}_k\|^2\}}.
	\]
	Therefore, the following bound holds:
	\begin{equation}
	\label{eq:sds:sumz1v1w}
	    \sum_{j=0}^{k-1} z_j (1-v_j(\kalpha)) (1-w_j)
	    \le \frac{\gammainc^4 \beta}{\mu\gammadec^4\min\{
		\gammainc^{-2}\alpha_0^2,
		\bar{\alpha}^2 \|\tilde{\vg}_k\|^2\}} 
		+ \frac{1}{\mu}\sum_{j=0}^{k-1} (1-z_j) w_j
		+ \log_{\gammadec^{-1}}\left(\frac{\alpha_0}{\gammadec\kalpha} \right).
	\end{equation}

    To conclude the proof, we combine~\eqref{eq:sds:sumzdecomp}, 
    \eqref{eq:sds:sumz1vw}, \eqref{eq:sds:sumzvw} and~\eqref{eq:sds:sumz1v1w} 
    as follows:
    \begin{eqnarray*}
    	\sum_{j=0}^{k-1} z_j 
    	&\le 
    	&\sum_{j=0}^{k-1} z_j (1-v_j(\kalpha))w_j
		+ \sum_{j=0}^{k-1} z_j v_j(\kalpha) w_j 
		+ \sum_{j=0}^{k-1} z_j (1-v_j(\kalpha)) (1-w_j) \\
		&\le
		&\frac{\gammainc^2 \beta}{\gammadec^2\min\{
		\gammainc^{-2}\alpha_0^2,
		\bar{\alpha}^2 [\mingk]^2\}}
		+\mu \sum_{j=0}^{k-1} (1-z_j)(1-w_j)		
		+\frac{\gammainc^4 \beta}{\mu\gammadec^4\min\{
		\gammainc^{-2}\alpha_0^2,
		\bar{\alpha}^2 [\mingk]^2\}} \\
		&
		&+ \frac{1}{\mu}\sum_{j=0}^{k-1} (1-z_j) w_j
		+ \log_{\gammadec^{-1}}\left(\frac{\alpha_0}{\gammadec\kalpha} \right) \\
		&=
		&\frac{(\mu\gammadec^2+\gammainc^2)\gammainc^2 \beta}{ \mu\gammadec^4\min\{
		\gammainc^{-2}\alpha_0^2,
		\bar{\alpha}^2 [\mingk]^2\}} 
		+ \mu \sum_{j=0}^{k-1} (1-z_j)(1-w_j) \\
		&
		&+ \frac{1}{\mu}\sum_{j=0}^{k-1} (1-z_j) w_j 
		+ \log_{\gammadec^{-1}}\left(\frac{\alpha_0}{\gammadec\kalpha} \right) \\
		&\le 
		&\frac{(\mu\gammadec^2+\gammainc^2)\gammainc^2 \beta}{ \mu\gammadec^4\min\{
		\gammainc^{-2}\alpha_0^2,
		\bar{\alpha}^2 [\mingk]^2\}} 
		+ \max\left\{\mu,\frac{1}{\mu}\right\} \sum_{j=0}^{k-1}(1-z_j) 
		+ \log_{\gammadec^{-1}}\left(\frac{\alpha_0}{\gammadec\kalpha} \right) \\
		&=
		&\frac{C}{\min\{\gammainc^{-2}\alpha_0^2,
		\bar{\alpha}^2 [\mingk]^2\}}
		+ \max\left\{\mu,\frac{1}{\mu}\right\} \sum_{j=0}^{k-1}(1-z_j)
		+ \log_{\gammadec^{-1}}\left(\frac{\alpha_0}{\gammadec\kalpha} \right),
    \end{eqnarray*}
    where the last line simply uses the formula~\eqref{eq:sds:constant} for
    $C$. Re-arranging the terms gives
    \[
    	\sum_{j=0}^{k-1} z_j 
    	\le 
    	\frac{\max\left\{\mu,\tfrac{1}{\mu}\right\}}
    	{1+\max\left\{\mu,\tfrac{1}{\mu}\right\}} k 
    	+\frac{1}{1+\max\left\{\mu,\tfrac{1}{\mu}\right\}}  
    	\left[\frac{C}{\min\{\gammainc^{-2}\alpha_0^2,
		\bar{\alpha}^2 [\mingk]^2\}}
		+ \log_{\gammadec^{-1}}\left(\frac{\alpha_0}{\gammadec\kalpha} \right)
		\right].
    \]
    By using the formula~\eqref{eq:sds:kalpha} for $\kalpha$ and observing that
    \begin{eqnarray*}
    	\frac{\max\left\{\mu,\tfrac{1}{\mu}\right\}}
    	{1+\max\left\{\mu,\tfrac{1}{\mu}\right\}}
    	&= &\max\left\{\frac{\mu}{1+\mu},\frac{1}{1+\mu}\right\}  
    	= p_0
    \end{eqnarray*}
    we arrive at the desired result.
\end{proof}

Note that Proposition~\ref{pr:sds:bndz} yields an alternate definition of 
$p_0$ than that identified by Gratton et 
al~\cite{SGratton_CWRoyer_LNVicente_ZZhang_2019}, which was 
$\tfrac{\log(\gammadec)}{\log(\gammainc^{-1}\gammadec)}$. However, we point out
that both terms in~\eqref{eq:sds:p0} are equal whenever 
$\gammainc = \gammadec^{-1}$, corresponding to $p_0=\tfrac{1}{2}$: this 
particular setting has been widely used in analyzing derivative-free 
methods based on probabilistic properties~\cite{CCartis_KScheinberg_2018}.

In addition to the result of Proposition~\ref{pr:sds:bndz}, we also 
have the following concentration inequality on the sum of the $z_k$ variables.

\begin{lemma}
\label{le:sds:conineq}
	Consider the sequences \revised{$\{\calD_k\}_k$} and $\{\mP_k\}_k$ generated by 
	Algorithm~\ref{algo:sds}, and suppose that the sequences are 
	$(\kappa,\dmax,p)$-descent and $(\eta,\sigma,\pmax,q)$-well-aligned, 
	respectively.
	
	Let $\pi_k(\lambda):=\P{\sum_{j=0}^{k-1} Z_j \le \lambda k}$ 
	for any $\lambda \in (0,pq)$. Then, 
	\begin{equation}
	\label{eq:sds!conineq}
		\pi_k(\lambda) \le \exp\left[-\frac{(pq-\lambda)^2}{2 pq}k\right].
	\end{equation}
\end{lemma}

The reasoning behind this result is the same \revised{as} Lemma 4.5 of Gratton 
et al.~\cite{SGratton_CWRoyer_LNVicente_ZZhang_2015}, except that the 
variable $Z_k$ involves two random events. To connect it with our 
probabilistic properties of interest, one must use
\begin{eqnarray*}
	\P{Z_k=1 | Z_0,\dots,Z_{k-1}} 
	&= 
	&\P{\calD_k\ \mbox{is ($\kappa,\dmax$)-descent}\ \&\ 
	\revised{\mP_k\ \mbox{is ($\eta,\sigma,\pmax$) well-aligned}}|\calF_{k-1}} \\
	&=
	&\P{\calD_k\ \mbox{is ($\kappa,\dmax$)-descent} |\calF_{k-1/2}} \\
	& 
	&\times \P{\revised{\mP_k\ \mbox{is ($\eta,\sigma,\pmax$) well-aligned}} |\calF_{k-1}},
\end{eqnarray*}
which holds because the events on $\calD_k$ and $\mP_k$ are conditionally 
independent.
%
%

Combining the results of Proposition~\ref{pr:sds:bndz} and 
Lemma~\ref{le:sds:conineq} 
lead to our main high-probability complexity result. The proof of this 
result follows from that of Theorem 4.6 in Gratton et 
al.~\cite{SGratton_CWRoyer_LNVicente_ZZhang_2015}, and merely differs 
by the presence of the additional logarithmic term in~\eqref{eq:sds:sumindic}.

\begin{theorem}
\label{th:highprobawcc}
	Suppose that the sequences \revised{$\{\calD_k\}_k$} and $\{\mP_k\}_k$ generated by 
	Algorithm~\ref{algo:sds} are $(\kappa,\dmax,p)$-descent and 
	$(\eta,\sigma,\pmax,q)$-well-aligned, respectively. Suppose further that 
	$pq>p_0$, where $p_0$ is defined as in Proposition~\ref{pr:sds:bndz}, and 
	consider an index $k$ and a tolerance $\eps>0$ such that
	\begin{equation}
	\label{eq:largekwcc}
		k \ge \frac{2}{pq-p_0}\left[ 
		\frac{(1-p_0)C}{\bar{\alpha}^2}\eps^{-2}
		+(1-p_0)\log_{\gammadec}\left(\tfrac{\gammadec^2 \bar{\alpha}}
		{\gammainc\alpha_0}\eps\right)
		\right]
	\end{equation}
	with $C$ defined as in Proposition~\ref{pr:sds:bndz}, and
	\begin{equation}
	\label{eq:smallepswcc}
		\eps \le \min\left\{1,
		\frac{\gammainc^{-1}\alpha_0}
		{\bar{\alpha}} \right\}.
	\end{equation}
	Then,
	\begin{equation}
	\label{eq:highprobawcc}
		\P{\tilde{G}_k \le \eps} \ge 1 - \exp\left[ 
		-\frac{(pq-p_0)^2}{8pq} k \right],
	\end{equation}
	where $\tilde{G}_k$ is the random variable associated with $\mingk$.
\end{theorem}

When $r=n$ and $\mP_k=\mI_n$, our result is of the same order \revised{as}
that of Gratton et al~\cite[Theorem 4.6]{SGratton_CWRoyer_LNVicente_ZZhang_2015}. 
The additional logarithmic term in~\eqref{eq:largekwcc} can be 
viewed as an additional cost incurred by the generalization of 
the reasoning to handle unbounded directions and random subspaces. 

The result of Theorem~\ref{th:highprobawcc} can be \revised{stated in a number 
of alternate ways}, depending on the quantity of interest (gradient norm, number of 
iterations) and the type of result that is sought (high 
probability, fixed probability guarantee, in expectation). We 
provide below two results that are of particular interest to 
our approach. The first one is a high-probability complexity bound 
on the number of function evaluations required to reach an 
approximate stationary point.

\begin{corollary}
\label{co:sds:wccevalsproba}
	Under the assumptions of Theorem~\ref{th:highprobawcc}, let 
	$N_{\epsilon}$ be the number of function evaluations required by 
	Algorithm~\ref{algo:sds} to reach a point $\vx_k$ such that 
	$\|\nabla f(\vx_k)\| \le \epsilon$, where $\epsilon>0$ 
	satisfies~\eqref{eq:smallepswcc}. Then,
	\begin{equation}
	\label{eq:sds:wccevalsproba}
		\P{N_{\epsilon} \le \left\lceil 
		\frac{2m}{pq-p_0}\phi(\epsilon)
		\right\rceil} \ge 
		1 - \exp\left[ 
		-\frac{(pq-p_0)}{4pq}\phi(\epsilon)  \right]
	\end{equation}
	where $\phi(\epsilon)= \frac{(1-p_0)C} {\bar{\alpha}^2}\eps^{-2}
	+(1-p_0)\log_{\gammadec}\left(\tfrac{\gammadec^2 \bar{\alpha}}
	{\gammainc\alpha_0}\eps\right)$.
\end{corollary}

The proof of Corollary~\ref{co:sds:wccevalsproba} follows from 
that of Theorem 4.8 in Gratton et 
al.~\cite{SGratton_CWRoyer_LNVicente_ZZhang_2015}. Its result shows that
Algorithm~\ref{algo:sds} has a high-probability complexity bound in
\begin{align}
    \frac{2m}{pq-p_0}\phi(\epsilon) = \calO\left(m \frac{1}{pq-p_0} 
	 \eta^{-2} \sigma^{-2} \pmax^4 \dmax^8 \kappa^{-2} \epsilon^{-2}\right) \label{eq_asymptotic_eval_bound}
\end{align}

The second corollary of our main result illustrates how our analysis leads to 
bounds in expectation: this result can be obtained following the argument 
developed for derivative-free algorithms based on probabilistic 
properties~\cite[Theorem 2.14]{SGratton_CWRoyer_LNVicente_ZZhang_2018}.

\begin{corollary}
\label{co:sds:wccevalsexp}
	Under the assumptions of Theorem~\ref{th:highprobawcc}, 
	\begin{equation}
	\label{eq:sds:wccevalsexp}
		\E{N_{\epsilon}} \le \frac{2 m}{pq-p_0} \phi(\epsilon) 
		+ \frac{1}{1-\exp\left(-\frac{(pq-p_0)^2}{8pq}\right)},
	\end{equation}
	where $N_{\epsilon}$ and $\phi(\epsilon)$ are defined as in 
	Corollary~\ref{co:sds:wccevalsproba}.
\end{corollary}

As $\phi(\epsilon) = \calO(\epsilon^{-2})$ for $\epsilon<1$, we obtain a 
complexity bound that matches that of other probabilistic techniques in terms 
of dependencies on $\epsilon$~\cite{SGratton_CWRoyer_LNVicente_ZZhang_2018}.
In addition, the dependencies on $m$, $\kappa$ and $\epsilon$ match that 
obtained for direct search based on deterministic 
descent~\eqref{eq:wccdsdeterm}. In the next two sections, we will establish 
evaluation complexity bounds for several choices of subspaces and polling 
sets.  

\subsection{Examples of direction generation techniques}
\label{subsec:sds:ps}

The above analysis of Algorithm~\ref{algo:sds} requires that our poll 
directions \revised{$\{\calD_k\}_k$} are $(\kappa,\dmax,p)$-descent and that our 
subspace matrices $\{\mP_k\}_k$ are $(\eta,\sigma,\pmax,q)$-well-aligned.
We now discuss several approaches to satisfying these requirements, both 
deterministic and probabilistic.

We begin by noting that our framework encompasses direct search based on 
deterministic and probabilistic descent. Indeed, if we take $r=n$ and 
$\mP_k=\mI_n$ for every $k$, $\mP_k$ is $(\eta,\sigma,\pmax,q)$-well-aligned with 
$\eta=\sigma=\pmax=q=1$, and Algorithm~\ref{algo:sds} reduces to 
Algorithm~\ref{algo:ds}. We then recover classical, deterministic direct search by choosing 
all \revised{$\calD_k$} to be the same PSS: the sequence \revised{$\{\calD_k\}_k$} is then 
$(\kappa,\dmax,p)$-descent with $p=1$, and $m \ge n+1$. Table 1 shows the 
values of $m$, $\kappa$ and $\dmax$ for three popular approaches: the 
coordinate vectors and their negatives, a set of $n+1$ vectors with uniform 
angles~\cite[Chapter 2.1]{ARConn_KScheinberg_LNVicente_2009b}, or the 
coordinate vectors and the vector with all entries $-1$.\footnote{We thank 
Warren Hare and Gabriel Jarry-Bolduc for checking the value of $\kappa$ in 
that latter case~\cite{WHare_GJarryBolduc_2020}.}
%
For probabilistic descent, we form \revised{$\calD_k$} by generating $m$ independent vectors 
uniformly distributed on the unit sphere: \revised{in} that case, \revised{$\calD_k$} is 
$(\kappa,\dmax,p)$-descent with $\dmax=1$, $\kappa=\tau/\sqrt{n}$ and 
$p \revised{\ge} 1-\left(\frac{1}{2}+\frac{\tau}{\sqrt{2\pi}}\right)^m$ for any 
$\tau \in [0,\sqrt{n}]$~\cite[Appendix B]{SGratton_CWRoyer_LNVicente_ZZhang_2015}.
Choosing $\tau=1$, we can make \revised{the lower bound on $p$ sufficiently large 
so that it exceeds $p_0$}
by taking $m=\bigO(1)$. These choices are summarized in 
Table~\ref{tab_pss_methods}.

\begin{table}[tb]
    \centering
    {\small
    \begin{tabular}{c|ccccc}
         \hline Method & $m$ & $\kappa$  & $\dmax$ & Success prob.~$p$ \\ \hline
         $[\mId_r,-\mId_r]$ & $2r$ & $r^{-1/2}$ & 1 & 1 \\
         Uniform angle PSS & $r+1$ & $r^{-1}$ & 1 & 1 \\
         $[\mId_r,-\ve_r]$ & $r+1$ & $(r^2+2(r-1)\sqrt{r})^{-1/2}$ & $\sqrt{r}$  & 1 \\
         Random unit vectors & $\bigO(1)$ & $r^{-1/2}$ & 1 
         & $1-(1/2 + 1/\sqrt{2\pi})^m$. \\
         \hline
    \end{tabular}
    } 
    \caption{Summary of methods for generating a direction set \revised{$\calD_k$} in 
    $\R^r$.}
    \label{tab_pss_methods}
\end{table}

We now explain how to generate $\mP_k$ matrices that are probabilistically 
well-aligned (the previous paragraph gives an example of a deterministically 
well-aligned matrix).
%
Consider first the requirements \eqref{eq:sds:Pk:grad} 
and~\eqref{eq:sds:Pk:pmax}. Three possible approaches can be used:
\begin{itemize}
    \item If $\mP_k$ has entries which are i.i.d.~$\mathcal{N}(0,1/r)$ and 
    $r = \Omega((1-\eta)^{-2}|\log (1-q)|)$, then $\mP_k$ 
    satisfies \eqref{eq:sds:Pk:grad} and \eqref{eq:sds:Pk:pmax}~\cite[Theorem 2.13]{
    SBoucheron_GLugosi_PMassart_2013} with 
    $\pmax = \Theta(\sqrt{n/r})$ with high probability~\cite[Corollary 3.11]{ASBandeira_RvHansel_2016}.
    \item If $\mP_k$ is an $s$-hashing matrix\footnote{i.e.~every column of 
    $\mP_k$ has exactly $s$ nonzero entries at randomly selected locations, 
    each taking value $\pm 1/\sqrt{s}$ with independent probability $1/2$.} 
    with $s=\Theta((1-\eta)^{-1}|\log (1-q)|)$ and 
    $r = \Omega((1-\eta)^{-2}|\log (1-q)|)$, then $\mP_k$ 
    satisfies \eqref{eq:sds:Pk:grad} and \eqref{eq:sds:Pk:pmax} for 
    $\pmax=\sqrt{n}$, where the value of $\pmax$ comes from 
    $\|\mP_k\|_2 \leq \|\mP_k\|_F = \sqrt{n}$~\cite{DKane_JNelson_2014}. 
    \item If $\mP_k = \sqrt{n/r} \mId_{r\times n}\mQ^T $ where $\mId_{r\times n}$ 
    denotes the first $r$ rows of the $n\times n$ identity matrix $\mI_n$, and
    $\mQ\in\R^{n\times n}$ is the orthogonal factor in the QR decomposition 
    $\mZ=\mQ\mR \in \R^{n\times n}$ of a matrix $\mZ$ with i.i.d. standard 
    normal entries such that the diagonal entries of $\mR$ are positive, then 
    $\|\mP_k\|_2 \leq \sqrt{n/r}$, so $\mP_k$ satisfies~\eqref{eq:sds:Pk:pmax}. 
    Moreover, for a given $r$ and $\eta$, $\mP_k$ satisfies~\eqref{eq:sds:Pk:grad} 
    with probability $q=1-I_{\eta r/d}(r/2, (n-r)/2)$, where $I_{p}(\alpha,\beta)$ is the 
    regularized incomplete Beta function 
    \cite[Lemma 1]{DKozak_SBecker_ADoostan_LTenorio_2021}.
    Although this does not give a closed form for a suitable choice of $r$, 
    numerical evidence suggests that $r$ may be chosen independently of 
    $n$~\cite[Figure 1]{DKozak_SBecker_ADoostan_LTenorio_2021}.
\end{itemize}

We now focus on condition \eqref{eq:sds:Pk:PTd}, that corresponds
to bounding the smallest singular value of $\mP_k$ away from zero. For 
Gaussian matrices with entries in $\mathcal{N}(0,1/r)$ as above, the 
estimate $\sigma_{\min}(\mP_k^T) \sim \sqrt{n/r}-1$ holds with high 
probability~\cite[Theorem 1.1]{MRudelson_RVershynin_2009}. In the case of 
orthogonal subsampling (the third case above), we automatically have 
$\sigma_{\min}(\mP_k^T)=\sqrt{n/r}$. We are unaware of theoretical 
results showing that $\sigma_{\min}(\mP_k) = \Theta(\sqrt{n/r})$ for hashing 
matrices, we present numerical evidence that this estimate holds when $n$ 
is sufficiently large compared to $r$ in Appendix~\ref{app_hashing_sing_val}.

%
%
%
%
%
%
%

\begin{table}[tb]
    \centering
    {\small
    \begin{tabular}{c|ccccc}
         \hline Method & $r$ & Success prob.~$q$ & $\pmax$ & $\sigma$ & Comments \\ \hline
         Identity & $n$ & 1 & 1 & 1 & --- \\
         Gaussian & $\bigO(1)$ & $1-e^{-\frac{(1-\eta)^2}{Ar}}$ & $\Theta(\sqrt{n/r})$ & $\Theta(\sqrt{n/r})$ & $\pmax$, $\sigma$ valid with high prob. \\
         $s$-Hashing & $\bigO(1)$ & $1-e^{-\frac{(1-\eta)^2}{Ar}}$ & $\sqrt{n}$ & $\Theta(\sqrt{n/r})$ & $s=\Theta(\frac{|\log(1-q)|}{1-\eta})$, $\sigma$ from App.~\ref{app_hashing_sing_val} \\
         Orthogonal & $\bigO(1)$ & $1-I_{\eta r/d}(\frac{r}{2},\frac{n-r}{2})$ &  $\sqrt{n/r}$ & $\sqrt{n/r}$ & $r$ estimated from \cite[Figure 1]{DKozak_SBecker_ADoostan_LTenorio_2021}. \\
         \hline
    \end{tabular}
    } 
    \caption{Summary of methods for generating $\mP_k$. The value $A$ is used to 
    represent a universal constant.}
    \label{tab_subspace_methods}
\end{table}

Table~\ref{tab_subspace_methods} summarizes our results for the three random choices of 
$\mP_k$ described above: for such choices, since $r$ can be chosen as small compared to 
the ambient dimension $n$, we may choose \revised{$\calD_k$} to be a standard PSS in $\R^r$, such as  
\revised{the columns of} $[\mId_r\ -\mId_r]$. As result, our poll step tests the points 
$\vx_k \pm \alpha_k \vp_{k,i}$, where $\vp_{k,i}$ is the $i$-th row of $\mP_k$. In that 
case, we obtain directions that are very similar to those used in Gratton et 
al.~\cite{SGratton_CWRoyer_LNVicente_ZZhang_2015}. In fact, using $r=1$ and orthogonal 
$\mP_k$ corresponds to using $\{\pm\sqrt{n}\ \vv\}$, where $\vv\in\R^n$ is a randomly 
drawn unit vector. Up to a scaling constant, this recovers the method described in 
Gratton et al.~\cite[p.~1535]{SGratton_CWRoyer_LNVicente_ZZhang_2015}. Our framework 
is however more general: for instance, using a Gaussian matrix $\mP_k$ with $r=1$ leads 
to the directions $\{\pm\vv\}$, where $\vv\in\R^n$ is a vector with i.i.d. standard 
Gaussian components.

\begin{remark}
Compared to traditional direct search, the use of random subspaces defined by 
$\mP_k$ increases the per-iteration linear algebra cost of Algorithm~\ref{algo:sds}:
however, this cost depends on the ensemble from which $\mP_k$ is generated.
For instance, a Gaussian $\mP_k$ would have a per-iteration cost of $\bigO(nr)$ to 
generate the matrix, and $\bigO(mnr)$ to construct each $\mP_k^T \vd$, noting that 
$m=\bigO(r)$ for most methods considered above. By contrast, an orthogonal $\mP_k$ 
requires $\bigO(nr^2)$ to generate $\mP_k$ via QR factorization, and again 
$\bigO(mnr)$ to construct $\mP_k^T \vd$.
Finally, using hashing to generate $\mP_k$ would take approximately 
$\bigO(ns)$ work, depending on the specific method used to select the nonzero 
locations, and only $\bigO(m\,\operatorname{nnz}(\mP_k))=\bigO(mns )$ to 
construct $\mP_k^T \vd$ using sparse multiplication. Note that this cost is lower 
when \revised{$\calD_k$ contains sparse vectors, such as the columns of 
$[\mId_r\ -\mId_r]$.}
\end{remark}


\subsection{Evaluation complexity results} 
\label{sec:sds:evalwcc}

As illustrated in Section~\ref{subsec:sds:ps}, the framework given by 
Algorithm~\ref{algo:sds} may be implemented in various ways. We now compare 
several instances of this method using the evaluation complexity 
bound~\eqref{eq_asymptotic_eval_bound} derived in Section~\ref{subsec:sds:wcc}. 
We pay particular attention to the dependency on $n$ that arises as we 
instantiate the method, since our reduced space approach aims at reducing this 
dependency.

Table~\ref{tab_algo_complexities} highlights the dependency on $n$ for several 
variants of the method. Note that for classical direct search ($\mP_k=\mId_n$, 
\revised{$\calD_k$} deterministic), there is a substantial difference in complexity 
depending on the choice of \revised{$\calD_k$}. On the contrary, using a randomized 
subspace choice $\mP_k$ always leads to an evaluation complexity in 
$\bigO(n\epsilon^{-2})$, independently of the choice of \revised{$\calD_k$}. 
Our result recovers the $\bigO(n\epsilon^{-2})$ complexity from direct search 
based on probabilistic descent~\cite{SGratton_CWRoyer_LNVicente_ZZhang_2015}, 
and shows how the same result may be achieved with substantially greater 
flexibility. It also gives further insight into the choice $\mP_k=\mId_n$ and 
$\revised{\calD_k}=\{\pm\vv\}$ where $\vv$ is a random unit vector, briefly discussed in 
Gratton et al.~\cite[p.~1535]{SGratton_CWRoyer_LNVicente_ZZhang_2015}. As 
explained above, this choice may be viewed as taking $\mP_k$ to be orthogonal 
with $r=1$ and \revised{$\calD_k$ as the columns of} $[\mId_n,-\mId_n]$.

\begin{table}[tb]
    \centering
    {\small
    \begin{tabular}{c|cccc}
         \hline 
         \revised{Columns of $\calD_k$} $\backslash$ $\mP_k$ choice 
         & Identity & Gaussian & $s$-Hashing & Orthogonal \\ 
         \hline
         $[\mId,-\mId]$ & $n^2$ & $n$ & $n$ & $n$ \\
         Uniform angle PSS & $n^3$ & $n$ & $n$ & $n$ \\
         $[\mId,-\ve]$ & $n^7$ & $n$ & $n$ & $n$ \\
         Random unit vectors & $n$ & $n$ & $n$ & $n$ \\
         \hline
    \end{tabular}
    } 
    \caption{Summary of evaluation complexity dependency on $n$ for 
    different choices of $\mP_k$ and \revised{$\calD_k$} in Algorithm~\ref{algo:sds}.}
    \label{tab_algo_complexities}
\end{table}


We conclude this section by commenting on the connections between our 
results and that obtained for model-based techniques. Unlike for direct search, 
in the model-based setting, using models that are of probabilistically 
good quality (or probabilistically fully 
linear~\cite{ARConn_KScheinberg_LNVicente_2009b}) is not known to improve the 
dependence on the dimension~\cite{SGratton_CWRoyer_LNVicente_ZZhang_2018}. 
However, such an improvement can be obtained by using deterministically fully 
linear models in randomly drawn subspaces~\cite{CCartis_LRoberts_2021}. This 
result, together with the complexity analysis of this section, suggests that 
using random subspace exploration is a more general approach method for 
producing the complexity improvement observed in both direct-search and 
model-based techniques.

\section{Numerical experiments}
\label{sec:num}

In this section, we investigate the practical performance of 
Algorithm~\ref{algo:sds} by comparing the following direct-search methods:
\begin{itemize}
    \item Algorithm~\ref{algo:ds} with deterministic descent, setting 
    $\calD_k$ as either the columns of $[\mId_n\ -\mId_n]$ or 
    $\calD_k=[\mId_n\ -\ve_n]$, where $\ve_n\in\R^n$ is the vector of ones;
    \item Algorithm~\ref{algo:ds} with probabilistic descent, using 
    $\calD_k = \{\pm \vv\}$ where $\vv\in\R^n$ is drawn from the 
    uniform distribution on the unit sphere $\|\vv\|=1$; 
    \item The Stochastic Three-Points (STP) algorithm from Bergou et 
    al.\cite{EBergou_EGorbunov_PRichtarik_2020}, that relies on opposite 
    directions uniform on the sphere, with step size schedule 
    $\alpha_k = \alpha_0 / (k+1)$;
    \item Algorithm~\ref{algo:sds} with $\mP_k$ taken either to be a 
    Gaussian matrix with $r\in\{1,2,3,4,5\}$, a hashing matrix, or a 
    subsampled orthogonal matrix with $r\in\{1,5\}$. The poll directions 
    were chosen to be columns of $[\mId_r\ -\mId_r]$.
\end{itemize}
All methods used $\alpha_0=1$ and terminated if $\alpha_k < 10^{-6}$.
The direct search methods used $\gammainc=2$, $\gammadec=0.5$, $\alpha_{\max}=1000$, and used opportunistic polling\footnote{i.e.~Do not check any other poll directions as soon as a $\vd_k$ satisfying \eqref{eq_new_sufficient_decrease} is found.} with the sufficient decrease condition
\begin{align}
    f(\vx_k + \alpha_k \mP_k^T \vd_k) < f(\vx_k) - \min\left(10^{-5}, 10^{-5} \alpha_k^2 \|\mP_k^T \vd_k\|^2\right), \label{eq_new_sufficient_decrease}
\end{align}
where $\mP_k=\mI$ for standard and probabilistic direct search.
All methods used in our experiments have been implemented within a 
Python package that has been made publicly 
available.\footnote{\url{https://github.com/lindonroberts/directsearch}}

\subsection{Robust Regression} 
\label{sec_robust_regression}

Our first set of experiments considers a robust linear regression problem with 
Gaussian data and a significant fraction of outliers, described by Carmon et 
al.~\cite{YCarmon_JDuchi_OHinder_ASidford_2017}. Specifically, we perform 
linear regression using a smoothed biweight loss function:
\begin{align}
    f(\vx) = \frac{1}{m}\sum_{i=1}^{m} \phi(\va_i^T\vx-b_i), 
    \qquad \text{where} \qquad \phi(\theta) := \frac{\theta^2}{1+\theta^2}. \label{eq_nonconvex_regression}
\end{align}
The vectors $\va_i\sim N(0,\mId_n)$ are i.i.d.~Gaussian vectors and the values 
$b_i$ are the coordinates of $\vb=\mA\vz + 3\vu_1+\vu_2$, where $\mA$ has 
columns $\va_i$, $\vz\sim N(0,4\mId_m)$, $\vu_1\sim N(0,\mId_m)$ and the entries of 
$\vu_2$ are i.i.d.~Bernoulli distributed with $p=0.3$. We choose $n=100$, 
$m=200$, start all solvers from the origin and run each solver for a maximum 
of $50(n+1)$ objective evaluations. For this initial study, we only compare 
direct search with $2n$ poll directions, probabilistic direct search, STP and 
Algorithm~\ref{algo:sds} with Gaussian $\mP_k$ and $r=1$. All randomized 
methods are run 10 times for a given problem. 

\begin{figure}[tbh]
	\centering
	\begin{subfigure}[b]{0.48\textwidth}
		\includegraphics[width=\textwidth]{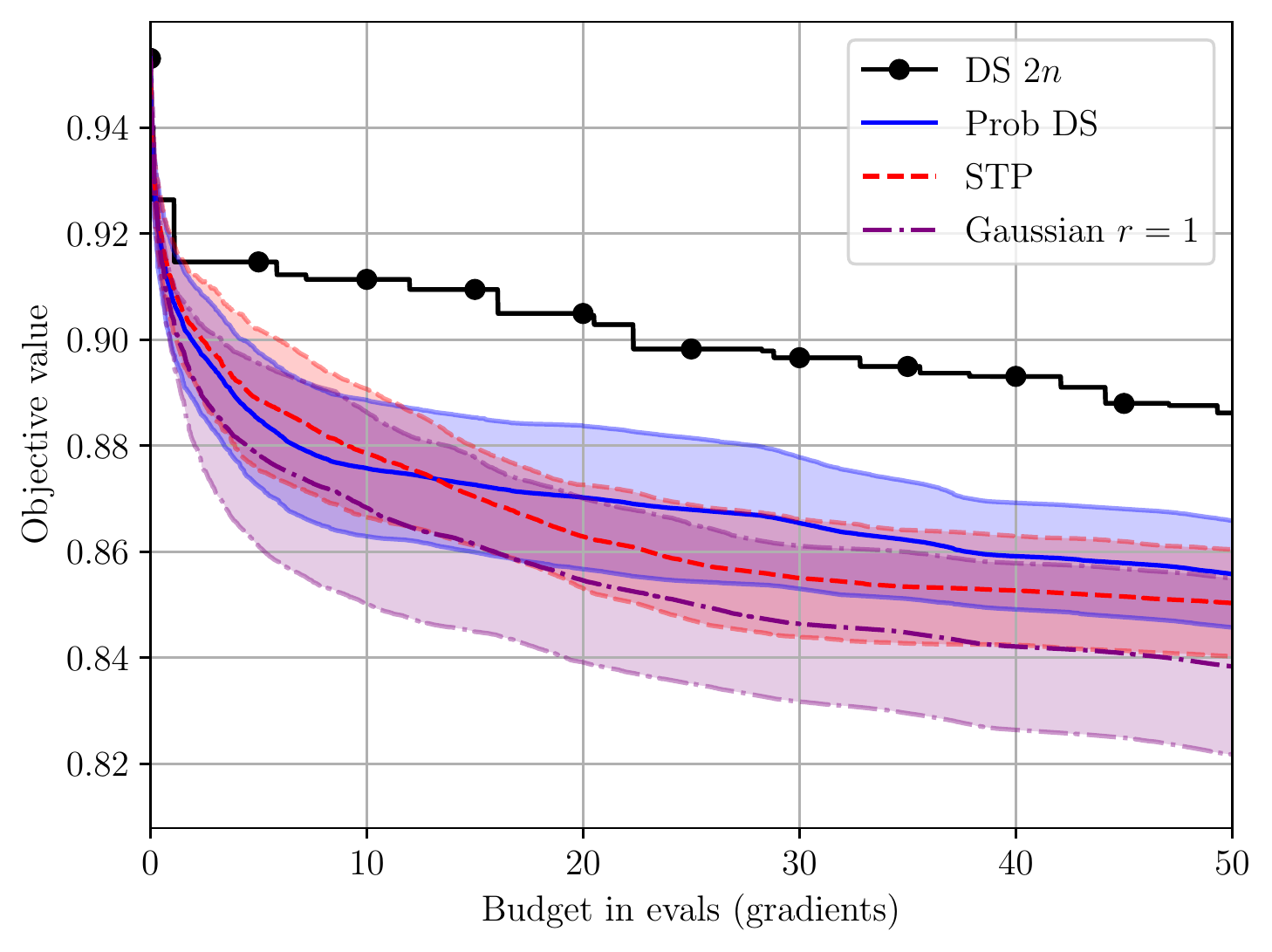}
		\caption{First problem instance}
	\end{subfigure}
	~
	\begin{subfigure}[b]{0.48\textwidth}
		\includegraphics[width=\textwidth]{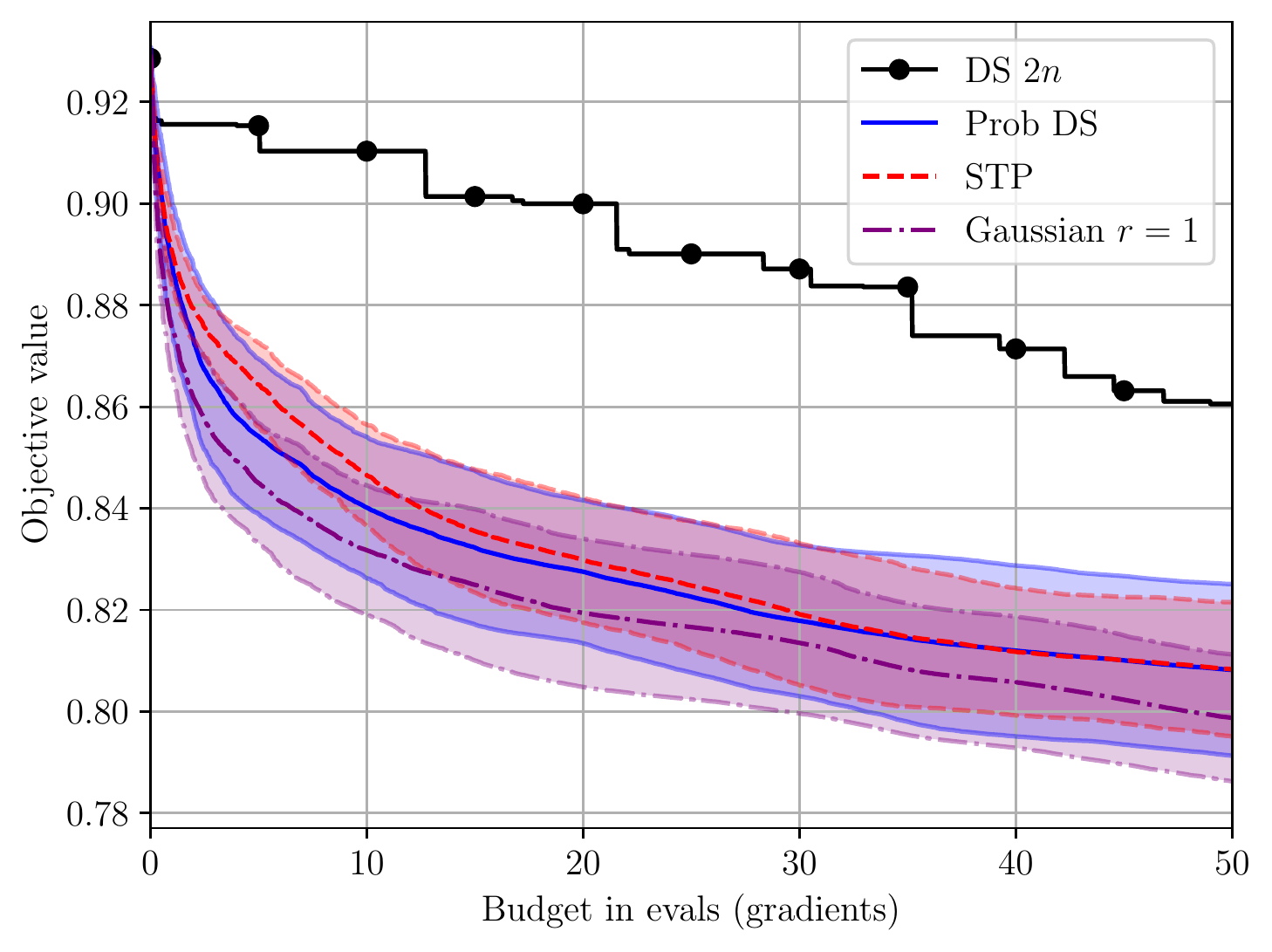}
		\caption{Second problem instance}
	\end{subfigure}
	\caption{Objective value versus number of objective evaluations (in units of $n+1$) for the nonconvex regression problem \eqref{eq_nonconvex_regression}. \revised{For randomized methods, we show the average value ($\pm 1$ standard deviation) achieved over 10 instances.}}
	\label{fig_nonconvex_regression}
\end{figure}

In Figure~\ref{fig_nonconvex_regression}, we plot the average objective 
value reached \revised{(over the 10 runs per solver with $\pm 1$ standard deviation error bounds)} 
versus the number of 
objective evaluations normalized in terms of simplex gradients (i.e.~units 
of $n+1$ evaluations). We show these results for two instances of 
problem~\eqref{eq_nonconvex_regression}, with different choices of $\va_i$ and 
$\vb$. All three randomized variants based on Algorithm~\ref{algo:sds} 
substantially outperform the traditional (deterministic) direct search method. 
Moreover, using randomized subspaces (effectively using Gaussian poll 
directions) leads to better results \revised{on average} than using unit vectors as poll 
directions, as done in probabilistic direct search and STP.

\subsection{CUTEst Collections} 
\label{sec_cutest}

Our second set of experiments is based on two collections of problems from 
the CUTEst optimization test 
set~\cite{NIMGould_DOrban_PhLToint_2015}:
\begin{itemize}
    \item The \textit{CFMR} collection \cite[Section 7]{CCartis_JFiala_BMarteau_LRoberts_2019} of 90 medium-scale problems ($25 \leq n \leq 120$, with $n\approx 100$ for 67 problems, approximately 75\%), ignoring bound constraints if any;
    \item The \textit{CR-large} collection \cite[Appendix B]{CCartis_LRoberts_2021} of 28 large-scale problems ($1000 \leq n \leq 5000$).
\end{itemize}
All solvers were run for at most $200(n+1)$ objective evaluations for the 
CFMR collection, and at most $10(n+1)$ evaluations for the CR-large collection.
As above, all randomized methods were run 10 times on every problem.

In this case, for each problem $P$ and solver instance $\calS$, we measure 
the number of objective evaluations $N_{P,\calS}$ required to achieve
\begin{align}
    f(\vx_k) \leq f^* + \tau(f(\vx_0) - f^*),
\end{align}
where $\tau\in(0,1)$ measures the required accuracy and $f^*$ is the true 
minimum objective value for each problem.
If this accuracy level is never achieved for a particular solver instance 
(within the desired budget or before termination), we take the convention 
$N_{P,\calS}=\infty$. We then compare solver performance using performance 
profiles \cite{EDolan_JMore_2002}, which, for a given solver $\calS$, plots 
the fraction of test problems for which $N_{P,\calS}$ is within some factor of 
$\min_{\calS} N_{P,\calS}$ (averaged over all solver instances).

\begin{figure}[t]
	\centering
	\begin{subfigure}[b]{0.48\textwidth}
		\includegraphics[width=\textwidth]{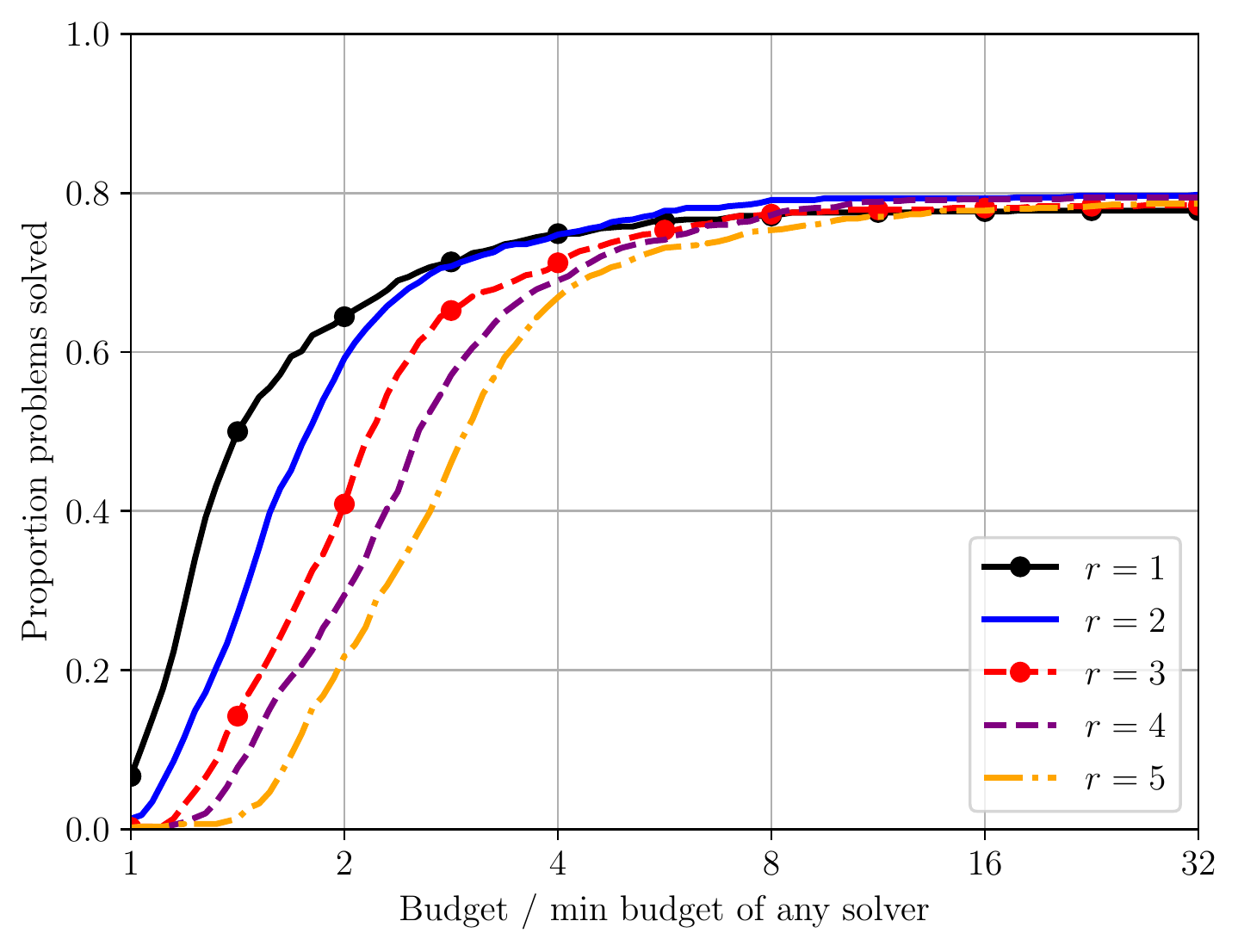}
		\caption{Medium-scale (CFMR), $\tau=10^{-1}$}
	\end{subfigure}
	~
	\begin{subfigure}[b]{0.48\textwidth}
		\includegraphics[width=\textwidth]{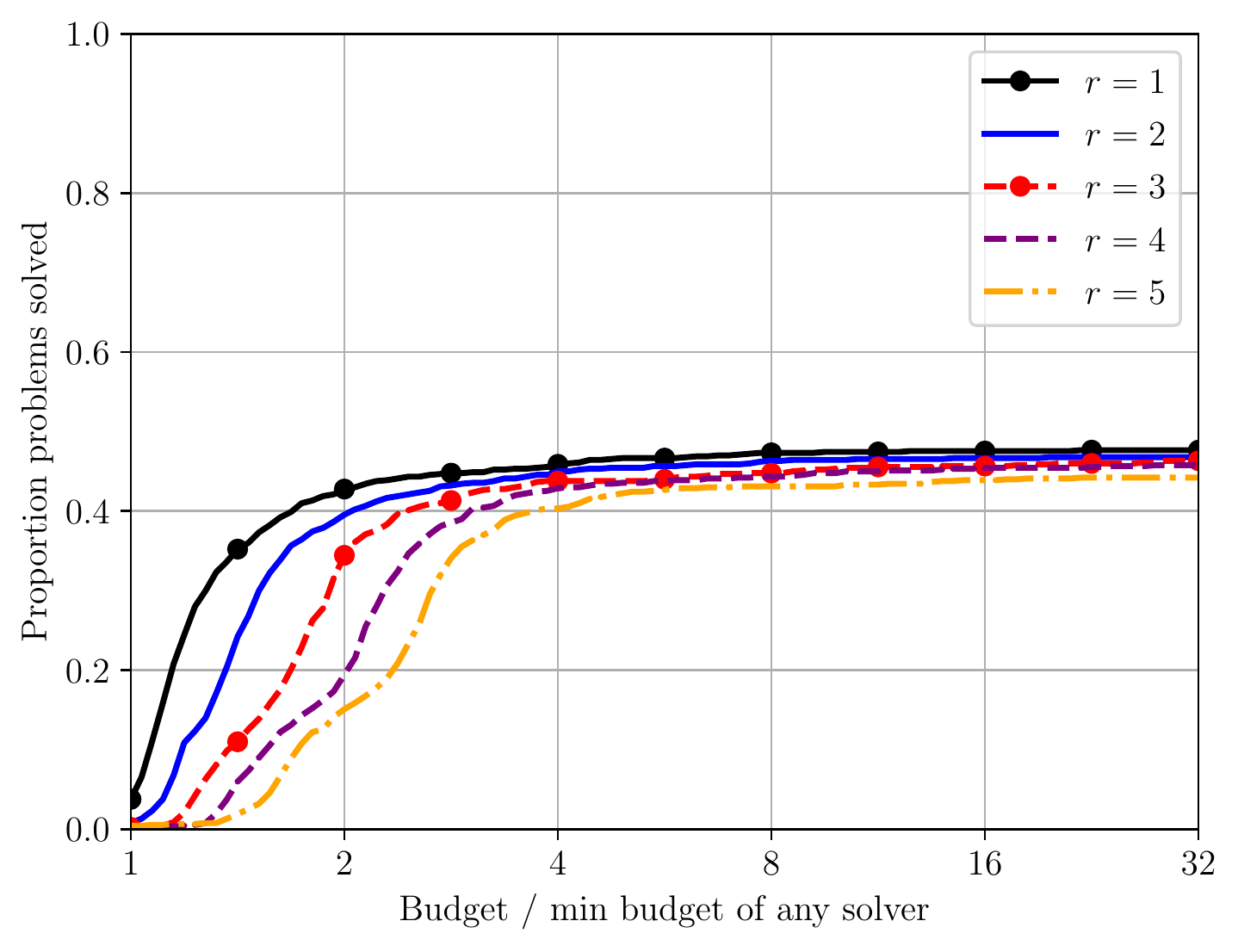}
		\caption{Medium-scale (CFMR), $\tau=10^{-3}$}
	\end{subfigure}
	\\
	\begin{subfigure}[b]{0.48\textwidth}
		\includegraphics[width=\textwidth]{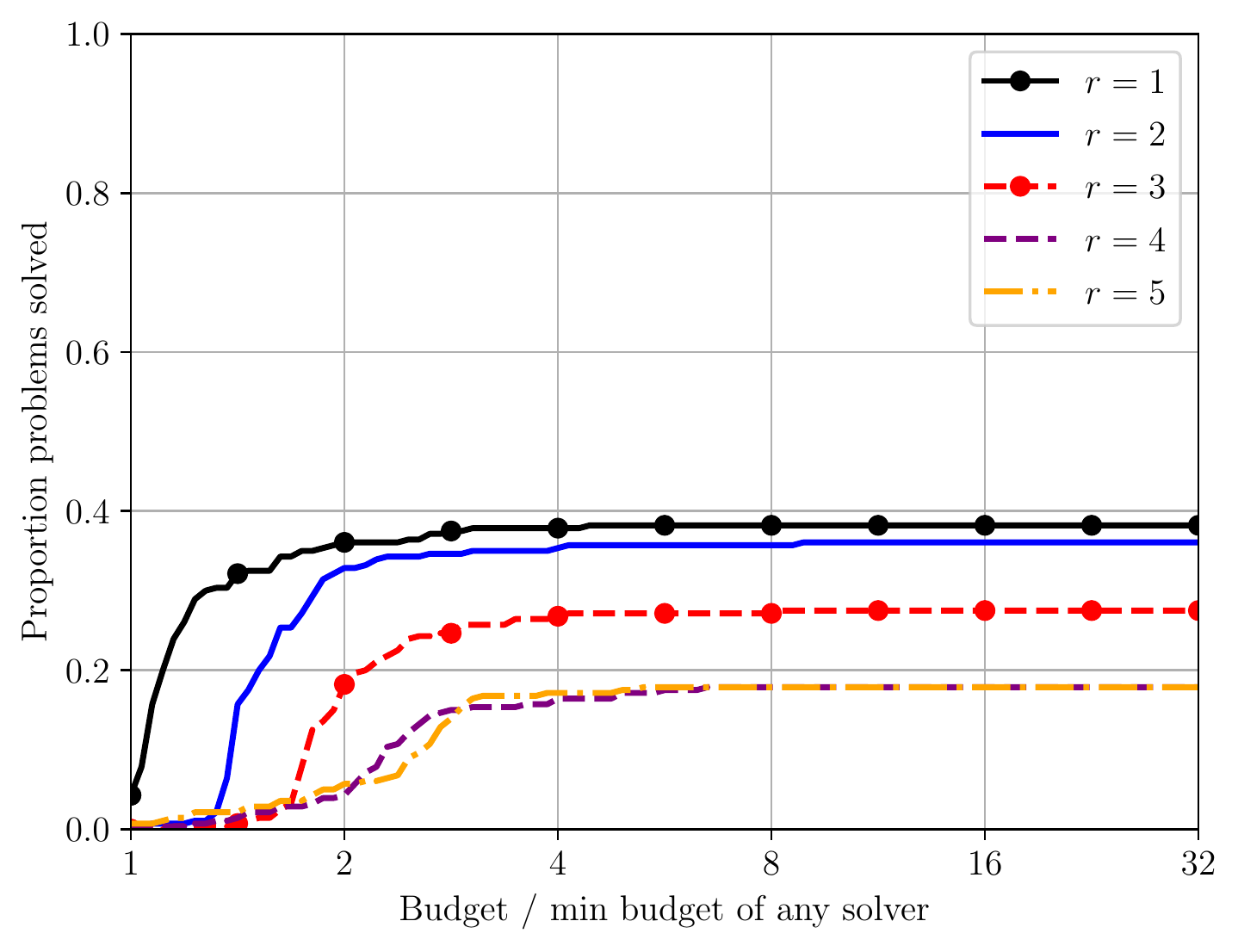}
		\caption{Large-scale (CR-large), $\tau=10^{-1}$}
	\end{subfigure}
	~
	\begin{subfigure}[b]{0.48\textwidth}
		\includegraphics[width=\textwidth]{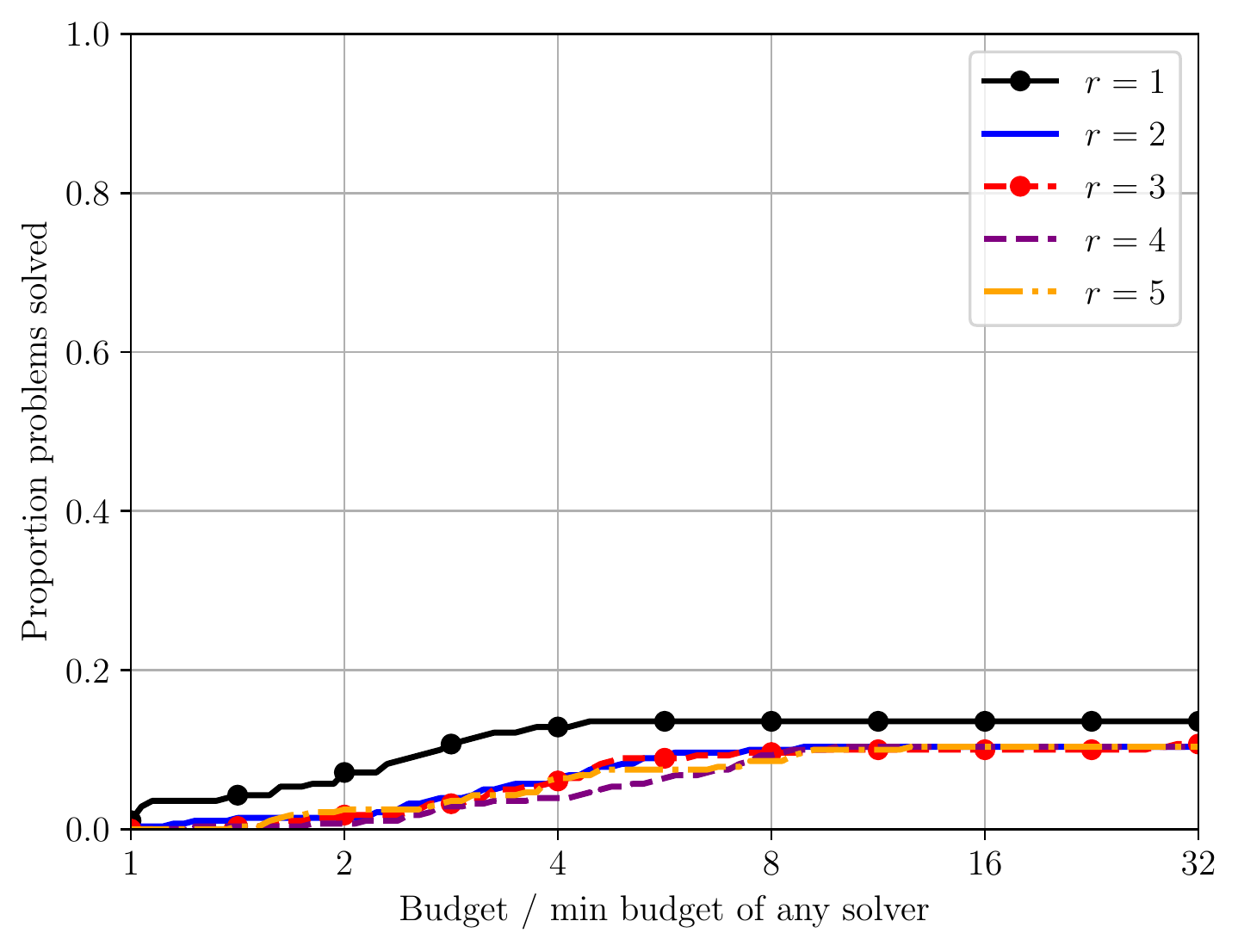}
		\caption{Large-scale (CR-large), $\tau=10^{-3}$}
	\end{subfigure}
	\caption{Comparison of $r$ values for Gaussian sketching}
	\label{fig_r_comparison}
\end{figure}

\begin{figure}[t]
	\centering
	\begin{subfigure}[b]{0.48\textwidth}
		\includegraphics[width=\textwidth]{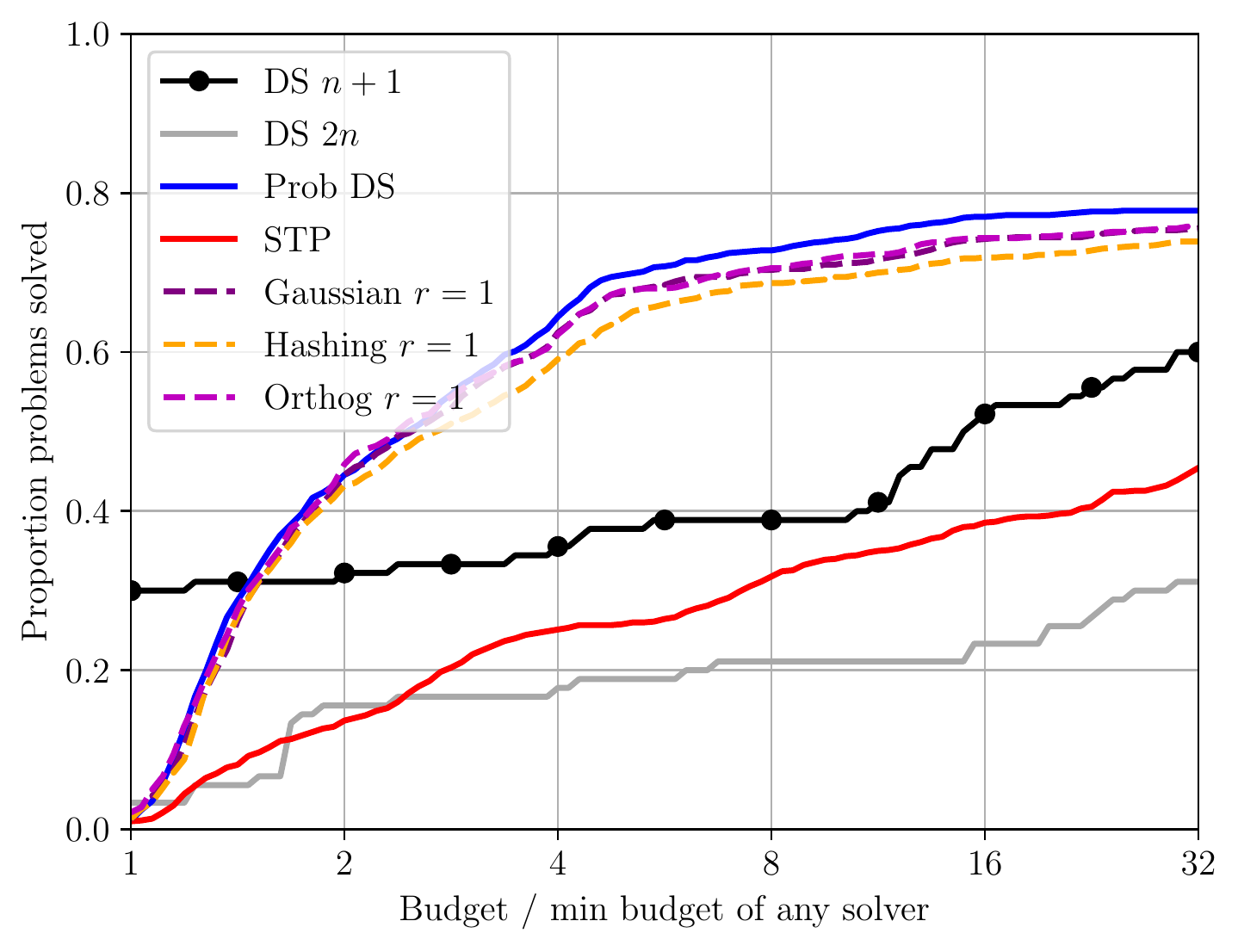}
		\caption{Medium-scale (CFMR), $\tau=10^{-1}$}
	\end{subfigure}
	~
	\begin{subfigure}[b]{0.48\textwidth}
		\includegraphics[width=\textwidth]{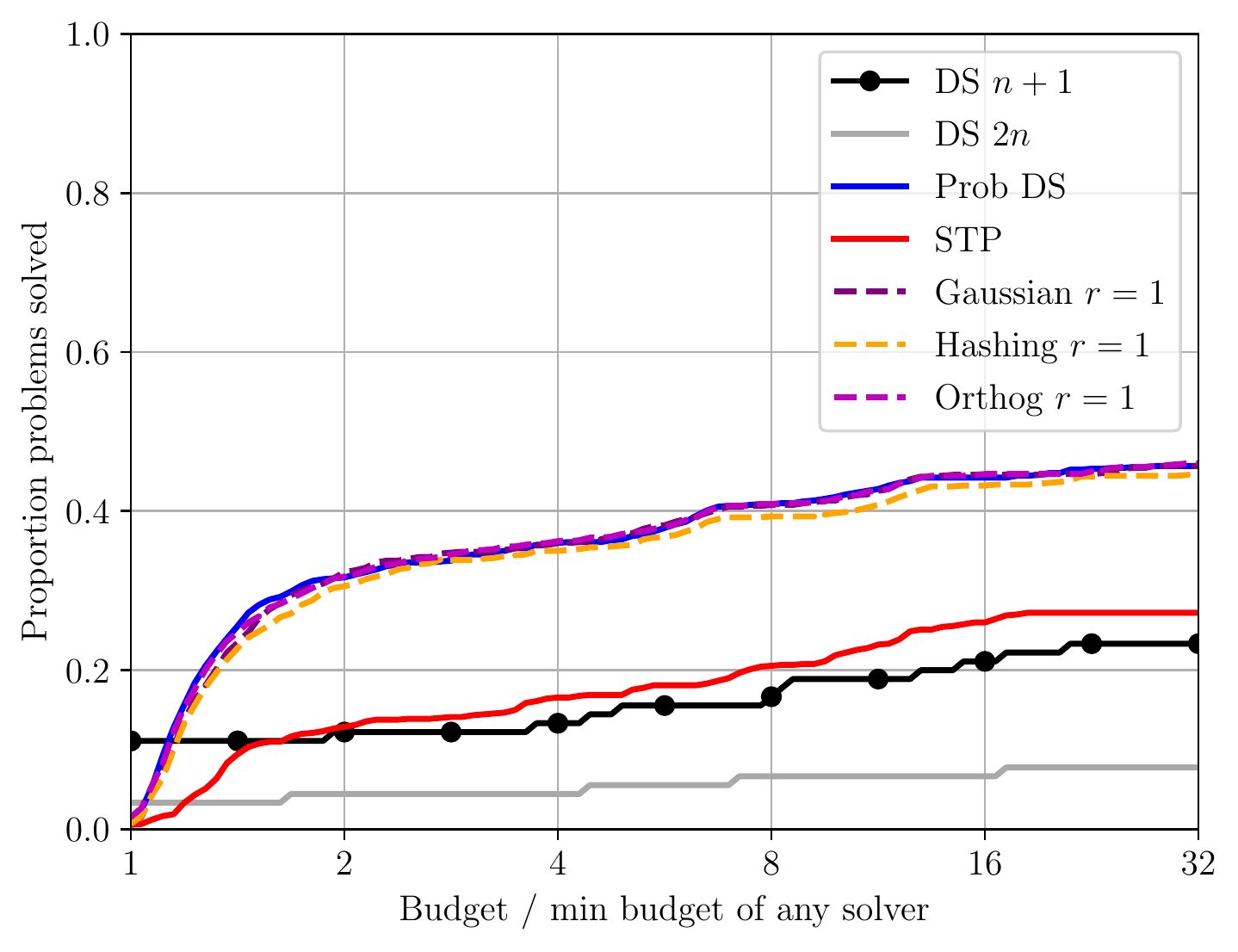}
		\caption{Medium-scale (CFMR), $\tau=10^{-3}$}
	\end{subfigure}
	\\
	\begin{subfigure}[b]{0.48\textwidth}
		\includegraphics[width=\textwidth]{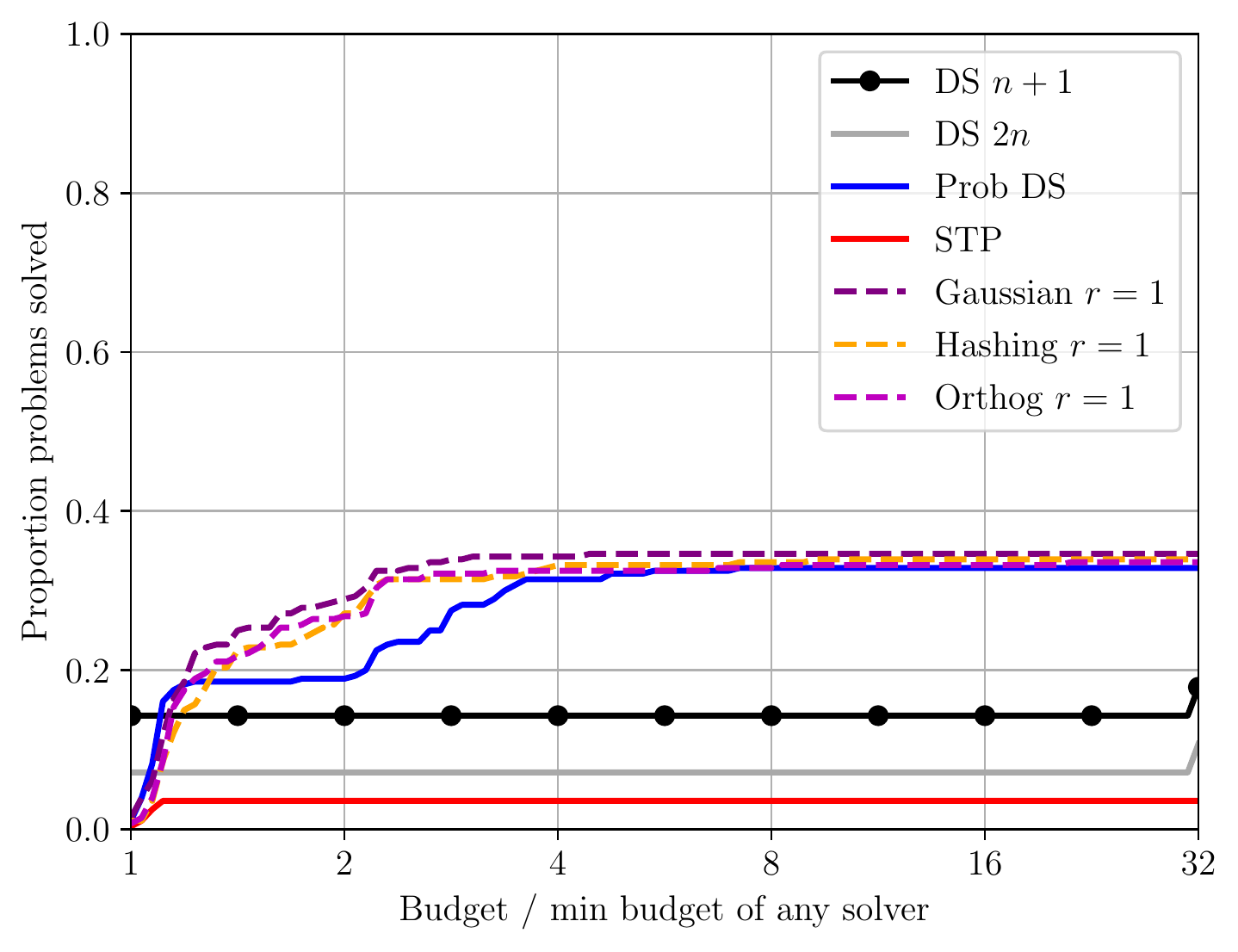}
		\caption{Large-scale (CR-large), $\tau=10^{-1}$}
	\end{subfigure}
	~
	\begin{subfigure}[b]{0.48\textwidth}
		\includegraphics[width=\textwidth]{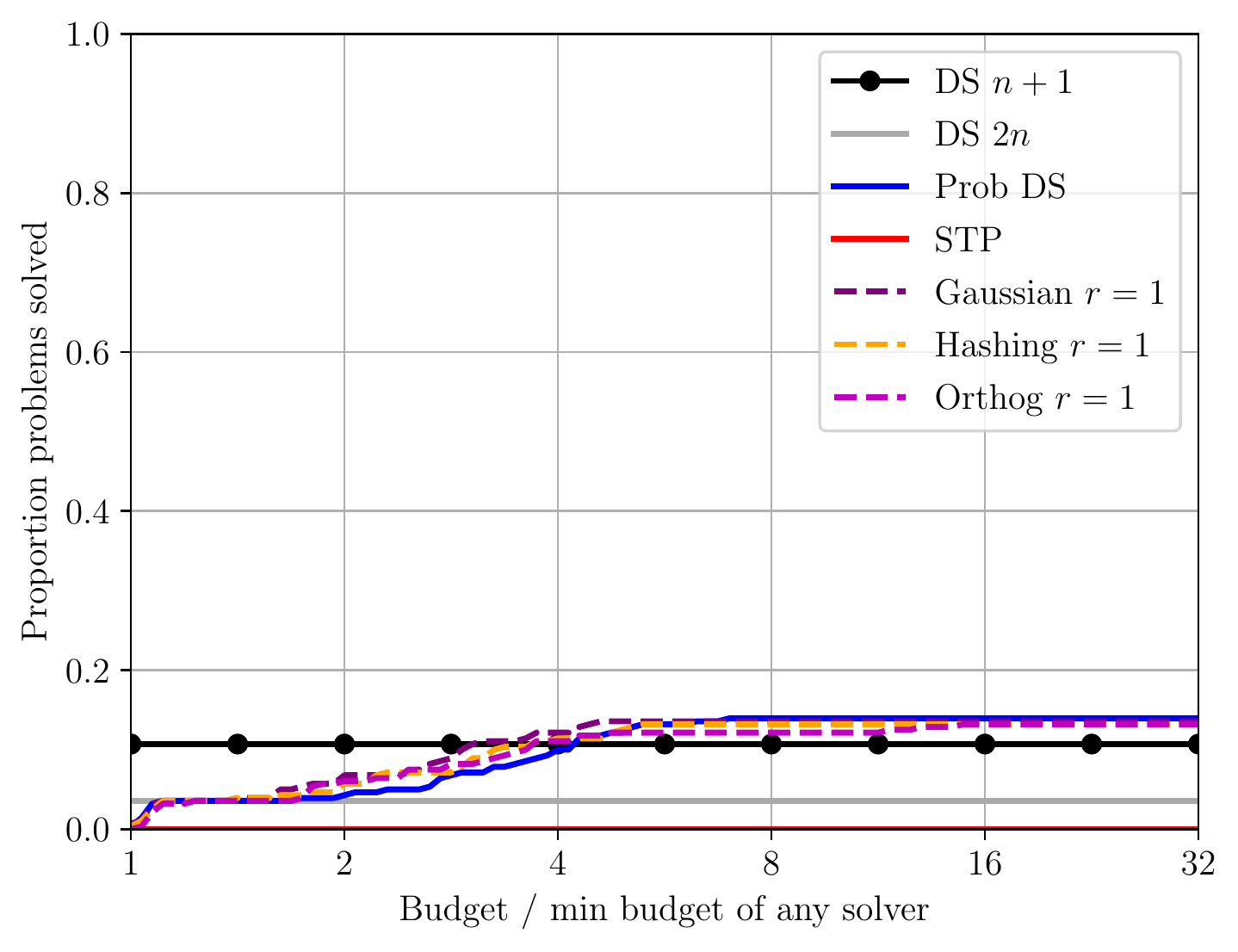}
		\caption{Large-scale (CR-large), $\tau=10^{-3}$}
	\end{subfigure}
	\caption{Comparison of different solvers}
	\label{fig_solver_comparison}
\end{figure}

Figure~\ref{fig_r_comparison} compares Algorithm~\ref{algo:sds} with Gaussian 
$\mP_k$ for different values of $r$, for both CFMR and CR-large test sets and 
accuracy levels $\tau\in\{10^{-1}, 10^{-3}\}$. In the case of the medium-scale 
CFMR problems, we see that smaller $r$ values give the best performance, but 
all values are able to solve essentially the same number of problems.
In the case of large-scale CR-large problems, again smaller $r$ values give 
the best performance, but we find that larger $r$ values are less robust in 
terms of the total number of problems solved (within the given evaluation 
budget).
Motivated by these results, in Figure~\ref{fig_solver_comparison} we compare 
Algorithm~\ref{algo:sds} (with Gaussian and hashing $\mP_k$ and $r=1$) with 
the other direct-search variants.
For both problem sets and accuracy levels, the best-performing methods are 
probabilistic direct search (Algorithm~\ref{algo:ds} with uniformly random 
directions) and the two Algorithm~\ref{algo:sds} variants, which for larger 
performance ratios strongly outperform deterministic direct search and STP.
Comparing Algorithm~\ref{algo:sds} with probabilistic direct search, we see 
broadly a very similar level of performance, reflecting the similarity of the 
two methods (both randomized with at most two objective evaluations per 
iteration at $\vx_k\pm\alpha_k \vv$ for a random vector $\vv$).
However, for low-accuracy solutions to the larger test problems CR-large, we 
see that allowing non-unit poll directions in Algorithm~\ref{algo:sds} gives 
a notable performance improvement over probabilistic direct search.

Overall, we find that Algorithm~\ref{algo:sds} performs similarly well to 
probabilistic direct search in many instances, although there appears to be a 
reasonable fraction of problems (in CUTEst and the nonconvex regression 
problem) for which our new approach gives superior performance.

\section{Conclusion}
\label{sec:conc}

We have proposed a general direct-search framework equipped with complexity 
guarantees that allows for polling in low-dimensional subspaces. We identified 
properties of both the subspaces and the polling directions that enable the 
derivation of complexity guarantees, which we expressed in a probabilistic 
form. Both our theoretical analysis and our numerical experiments suggest that 
using one-dimensional subspaces is an overall efficient approach \revised{that} 
allows for applying direct-search methods to larger dimensional problems than 
usually considered in the literature. Interestingly, the use of random embeddings 
does not remove all dependencies on the original dimension of the problem, 
but comes with the computational advantage of drawing only two directions 
per iteration.

Our approach is based on characterizing properties of the polling directions 
at the iteration level, while information from previous iterations is 
typically used in practice to improve the polling set. On the other hand, removing our 
conditional independence yields a number of mathematical challenges, that 
would have to be overcome in order to obtain a complexity result. Similarly, 
extending this method to handle noisy function values would require combining our analysis with existing frameworks for stochastic optimization\revised{, some of which have strong similarity to direct search~\cite{JCDuchi_MIJordan_MJWainwright_AWibisono_2015,JLarson_MMenickelly_SMWild_2019}.}
Both avenues of research represent natural perspectives for future work.

\bibliographystyle{plain} 
\bibliography{refsSDFO}

\begin{thebibliography}{10}

\bibitem{CAudet_WHare_2017}
C.~Audet and W.~Hare.
\newblock {\em Derivative-{F}ree and {B}lackbox {O}ptimization}.
\newblock Springer Series in Operations Research and Financial Engineering.
  Springer International Publishing, 2017.

\bibitem{ASBandeira_KScheinberg_LNVicente_2014}
A.~S. {Bandeira}, K.~Scheinberg, and L.~N. Vicente.
\newblock Convergence of trust-region methods based on probabilistic models.
\newblock {\em SIAM J. Optim.}, 24:1238--1264, 2014.

\bibitem{ASBandeira_RvHansel_2016}
A.~S. Bandeira and R.~{van Handel}.
\newblock Sharp nonasymptotic bounds on the norm of random matrices with
  independent entries.
\newblock {\em The Annals of Probability}, 44(4), 2016.

\bibitem{EBergou_EGorbunov_PRichtarik_2020}
E.~Bergou, E.~Gorbunov, and P.~Richt\'{a}rik.
\newblock Stochastic three points method for unconstrained smooth minimization.
\newblock {\em SIAM J. Optim.}, 30:2726--2749, 2020.

\bibitem{SBoucheron_GLugosi_PMassart_2013}
S.~Boucheron, G.~Lugosi, and P.~Massart.
\newblock {\em Concentration Inequalities: A Nonasymptotic Theory of
  Independence}.
\newblock Oxford University Press, Oxford, 2013.

\bibitem{YCarmon_JDuchi_OHinder_ASidford_2017}
Y.~Carmon, J.~C Duchi, O.~Hinder, and A.~Sidford.
\newblock ``{C}onvex until proven guilty'': Dimension-free acceleration of
  gradient descent on non-convex functions.
\newblock In {\em Proceedings of the 34th International Conference on Machine
  Learning}, volume~70, Sydney, 2017. PMLR.

\bibitem{CCartis_TFerguson_LRoberts_2020}
C.~Cartis, T.~Ferguson, and L.~Roberts.
\newblock Scalable derivative-free optimization for nonlinear least-squares
  problems.
\newblock In {\em ICML Workshop on Beyond First Order Methods in ML Systems},
  2020.

\bibitem{CCartis_JFiala_BMarteau_LRoberts_2019}
C.~Cartis, J.~Fiala, B.~Marteau, and L.~Roberts.
\newblock Improving the flexibility and robustness of model-based
  derivative-free optimization solvers.
\newblock {\em ACM Trans. Math. Software}, 45:32:1--32:41, 2019.

\bibitem{CCartis_JFowkes_ZShao_2020}
C.~Cartis, J.~Fowkes, and Z.~Shao.
\newblock A randomised subspace {G}auss-{N}ewton method for nonlinear
  least-squares.
\newblock In {\em Workshop on ``Beyond first-order methods in ML systems'' at
  the 37th International Conference on Machine Learning}, Vienna, Austria,
  2020.

\bibitem{CCartis_JFowkes_ZShao_2022}
C.~Cartis, J.~Fowkes, and Z.~Shao.
\newblock Randomised subspace methods for non-convex optimization, with
  applications to nonlinear least-squares.
\newblock Technical report, University of Oxford, 2022.

\bibitem{CCartis_NIMGould_PhLToint_2012b}
C.~Cartis, N.~I.~M. {Gould}, and Ph.~L. {Toint}.
\newblock On the oracle complexity of first-order and derivative-free
  algorithms for smooth nonconvex minimization.
\newblock {\em SIAM J. Optim.}, 22:66--86, 2012.

\bibitem{CCartis_LRoberts_2021}
C.~Cartis and L.~Roberts.
\newblock Scalable subspace methods for derivative-free nonlinear least-squares
  optimization.
\newblock arXiv:2102.12016, 2021.

\bibitem{CCartis_KScheinberg_2018}
C.~Cartis and K.~Scheinberg.
\newblock Global convergence rate analysis of unconstrained optimization
  methods based on probabilistic models.
\newblock {\em Math. Program.}, 169:337--375, 2018.

\bibitem{ARConn_KScheinberg_LNVicente_2009b}
A.~R. {Conn}, K.~Scheinberg, and L.~N. {Vicente}.
\newblock {\em Introduction to Derivative-Free Optimization}.
\newblock MPS-SIAM Series on Optimization. Society for Industrial and Applied
  Mathematics, Philadelphia, 2009.

\bibitem{MDodangeh_LNVicente_ZZhang_2016}
M.~Dodangeh, L.~N. {Vicente}, and Z.~Zhang.
\newblock On the optimal order of worst case complexity of direct search.
\newblock {\em Optim. Lett.}, 10:699--708, 2016.

\bibitem{EDolan_JMore_2002}
E.~D. Dolan and J.~J. Mor{\'e}.
\newblock Benchmarking optimization software with performance profiles.
\newblock {\em Math. Program.}, 91(2):201--213, 2002.

\bibitem{JCDuchi_MIJordan_MJWainwright_AWibisono_2015}
J.~C. {Duchi}, M.~I. {Jordan}, M.~J. {Wainwright}, and A.~Wibisono.
\newblock Optimal rates for zero-order convex optimization: the power of two
  function evaluations.
\newblock {\em IEEE Trans. Inform. Theory}, 61:2788--2806, 2015.

\bibitem{RGarmanjani_DJudice_LNVicente_2016}
R.~Garmanjani, D.~J\'{u}dice, and L.~N. {Vicente}.
\newblock Trust-region methods without using derivatives: {W}orst case
  complexity and the non-smooth case.
\newblock {\em SIAM J. Optim.}, 26:1987--2011, 2016.

\bibitem{NIMGould_DOrban_PhLToint_2015}
N.~I.~M. {Gould}, D.~Orban, and Ph.~L. {Toint}.
\newblock {CUTE}st: a constrained and unconstrained testing environment with
  safe threads.
\newblock {\em Comput. Optim. Appl.}, 60:545--557, 2015.

\bibitem{SGratton_CWRoyer_LNVicente_ZZhang_2015}
S.~Gratton, C.~W. {Royer}, L.~N. {Vicente}, and Z.~Zhang.
\newblock Direct search based on probabilistic descent.
\newblock {\em SIAM J. Optim.}, 25:1515--1541, 2015.

\bibitem{SGratton_CWRoyer_LNVicente_ZZhang_2018}
S.~Gratton, C.~W. {Royer}, L.~N. {Vicente}, and Z.~Zhang.
\newblock Complexity and global rates of trust-region methods based on
  probabilistic models.
\newblock {\em IMA J. Numer. Anal.}, 38:1579--1597, 2018.

\bibitem{SGratton_CWRoyer_LNVicente_ZZhang_2019}
S.~Gratton, C.~W. {Royer}, L.~N. {Vicente}, and Z.~Zhang.
\newblock Direct search based on probabilistic feasible descent for bound and
  linearly constrained problems.
\newblock {\em Comput. Optim. Appl.}, 72:525--559, 2019.

\bibitem{WHare_GJarryBolduc_2020}
W.~Hare and G.~{Jarry-Bolduc}.
\newblock A deterministic algorithm to compute the cosine measure of a finite
  positive spanning set.
\newblock {\em Optim. Lett.}, 14:1305--1316, 2020.

\bibitem{DKane_JNelson_2014}
Daniel~M. Kane and Jelani Nelson.
\newblock Sparser {J}ohnson-{L}indenstrauss transforms.
\newblock {\em Journal of the ACM}, 61(1):1--23, 2014.

\bibitem{TGKolda_RMLewis_VTorczon_2003}
T.~G. {Kolda}, R.~M. {Lewis}, and V.~Torczon.
\newblock Optimization by direct search: {N}ew perspectives on some classical
  and modern methods.
\newblock {\em SIAM Rev.}, 45:385--482, 2003.

\bibitem{DKozak_SBecker_ADoostan_LTenorio_2021}
D.~Kozak, S.~Becker, A.~Doostan, and L.~Tenorio.
\newblock A stochastic subspace approach to gradient-free optimization in high
  dimensions.
\newblock {\em Comput. Optim. Appl.}, 79:339--368, 2021.

\bibitem{DKozak_CMolinari_LRosasco_LTenorio_SVilla_2021}
D.~Kozak, C.~Molinari, L.~Rosasco, L.~Tenorio, and S.~Villa.
\newblock Zeroth order optimization with orthogonal random directions.
\newblock arXiv:2107.03941v1, 2021.

\bibitem{JLarson_MMenickelly_SMWild_2019}
J.~Larson, M.~Menickelly, and S.~M. {Wild}.
\newblock Derivative-free optimization methods.
\newblock {\em Acta Numer.}, 28:287--404, 2019.

\bibitem{YuNesterov_2011}
Yu. Nesterov.
\newblock Random gradient-free minimization of convex functions.
\newblock Technical Report 2011/1, {CORE}, Universit\'{e} {C}atholique de
  {L}ouvain, 2011.

\bibitem{YuNesterov_VSpokoiny_2017}
Yu. Nesterov and V.~Spokoiny.
\newblock Random gradient-free minimization of convex functions.
\newblock {\em Found. Comput. Math.}, 17:527--566, 2017.

\bibitem{CPaquette_KScheinberg_2020}
C.~Paquette and K.~Scheinberg.
\newblock A stochastic line search method with convergence rate analysis.
\newblock {\em SIAM J. Optim.}, 30:349--376, 2020.

\bibitem{MRudelson_RVershynin_2009}
M.~Rudelson and R.~Vershynin.
\newblock Smallest singular value of a random rectangular matrix.
\newblock {\em Communications on Pure and Applied Mathematics},
  62(12):1707--1739, 2009.

\bibitem{ZShao_2022}
Z.~Shao.
\newblock {\em On Random Embeddings and their Applications to Optimization}.
\newblock PhD thesis, University of Oxford, 2022.

\bibitem{LNVicente_2013}
L.~N. {Vicente}.
\newblock Worst case complexity of direct search.
\newblock {\em EURO J. Comput. Optim.}, 1:143--153, 2013.

\end{thebibliography}

\appendix

\section{Smallest Singular Value of Hashing Matrices} 
\label{app_hashing_sing_val}

This appendix provides numerical evidence for the bound 
$\sigma = \Theta(\sqrt{n/r})$ when $\mP_k$ is a hashing matrix 
(see the description in Section~\ref{subsec:sds:ps}) and $\sigma>0$ is the 
parameter used in Definitions~\ref{de:sds:Pk} and~\ref{de:sds:pbaPk}.
To this end, we generated 100 independent hashing matrices with $s=1$ for 
different combinations of $n$ and $r$ (with $n\gg r$)  and calculated their 
smallest singular value. Note that for hashing matrices, we have 
$\sigma \leq \sigma_{\min}(\mP_k)=\sigma_{\min}(\mP_k^T) \leq \pmax$, and 
that $\pmax\leq \sqrt{n}$.

\begin{figure}[tb]
	\centering
	\begin{subfigure}[b]{0.48\textwidth}
		\includegraphics[width=\textwidth]{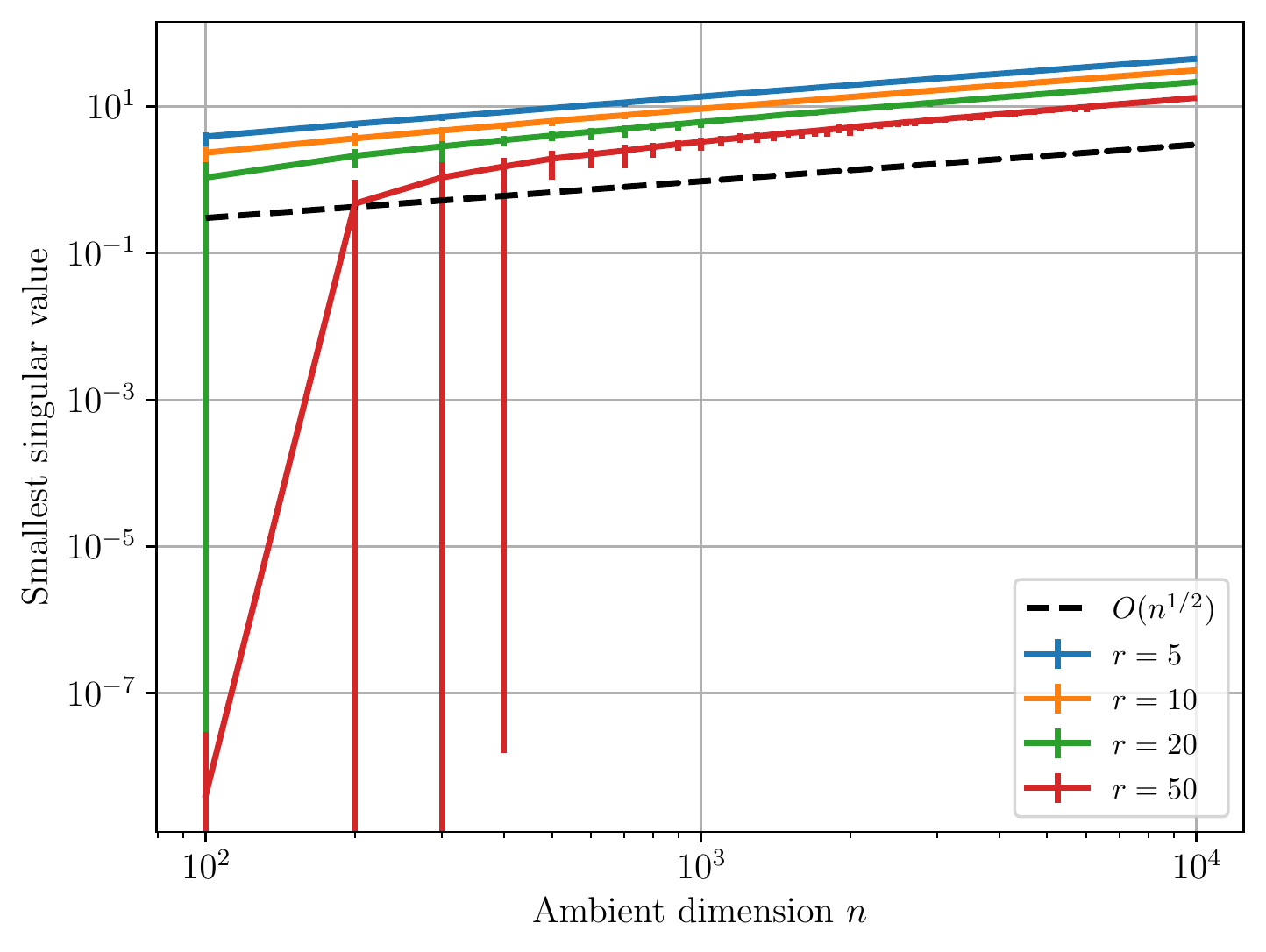}
		\caption{$\sigma_{\min}(\mP_k)$ as a function of $n$.}
		\label{fig_hashing_sing_vals_n}
	\end{subfigure}
	~
	\begin{subfigure}[b]{0.48\textwidth}
		\includegraphics[width=\textwidth]{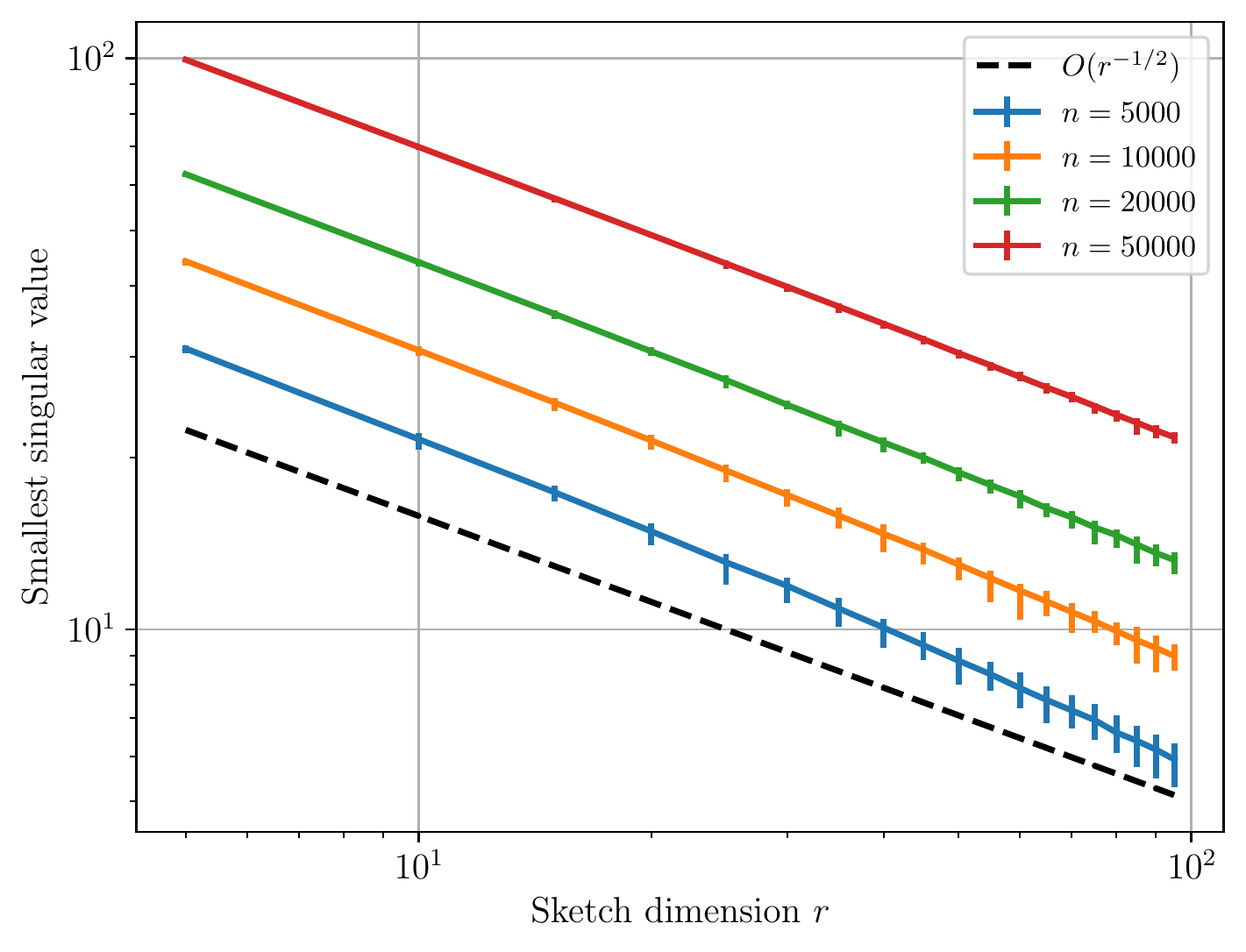}
		\caption{$\sigma_{\min}(\mP_k)$ as a function of $r$.}
		\label{fig_hashing_sing_vals_r}
	\end{subfigure}
	\caption{Computed values of $\sigma_{\min}(\mP_k)$ (minimum nonzero 
	singular value) for hashing with $s=1$, as a function of $n$ and $r$.}
	\label{fig_hashing_sing_vals}
\end{figure}

In Figure~\ref{fig_hashing_sing_vals}, we plot the mean, maximum and 
minimum values for each of the 100 independent trials. When $n \gg r$, we 
observe that there is very little uncertainty in $\sigma_{\min}(\mP_k)$, 
and the relationship $\sigma_{\min}(\mP_k) = \bigO(\sqrt{n})$ holds, as 
shown by Figure~\ref{fig_hashing_sing_vals_n}. Similarly, we see little 
uncertainty in $\sigma_{\min}(\mP_k)$ as $r$ varies (see 
Figure~\ref{fig_hashing_sing_vals_r}), while the relationship 
$\sigma_{\min}(\mP_k) = \bigO(1/\sqrt{r})$ holds. 
Altogether, these results suggest that the value 
$\sigma=\Theta(\sqrt{n/r})$ leads to satisfaction of~\eqref{eq:sds:Pk:PTd} 
with high probability when $\mP_k$ is a hashing matrix.

\end{document}